\documentclass[10pt,a4paper]{amsart}
\usepackage{setspace}
\usepackage{amsfonts}

\usepackage{amsthm}
\usepackage{amsmath}
\usepackage{amssymb}
\usepackage{amscd}
\usepackage[utf8]{inputenc}
\usepackage{t1enc}
\usepackage[mathscr]{eucal}
\usepackage{indentfirst}
\usepackage{graphicx}
\usepackage{graphics}
\usepackage{pict2e}
\usepackage{epic}
\usepackage{float}
\usepackage{booktabs}
\usepackage{MnSymbol}
\usepackage{multirow}
\usepackage{shuffle} %the symbol shuffle product
\usepackage{todonotes}
\numberwithin{equation}{section}
\usepackage[margin=2.9cm]{geometry}
\usepackage{epstopdf}
\usepackage[hooks]{tcolorbox}

\usepackage{tikz-cd}
\usepackage{quiver}

\makeatletter
\def\tcb@restore@footnote{%
	\def\@mpfn{footnote}%
	\def\thempfn{\arabic{footnote}}%
	\let\@footnotetext\tcb@footnote@collect
}

\long\def\tcb@footnote@collect#1{%
	\expandafter\gappto\expandafter\tcb@footnote@acc\expandafter{%
		\expandafter\footnotetext\expandafter[\@thefnmark]{#1}%
	}%
}

\def\tcb@footnote@use{%
	\tcb@footnote@acc
	\global\let\tcb@footnote@acc\@empty
}
\global\let\tcb@footnote@acc\@empty

\tcbset{
	every box/.style={
		before upper pre=\tcb@restore@footnote
	},
	every box on layer 1/.append style={
		after app=\tcb@footnote@use
	}
}
\makeatother

\theoremstyle{plain}
\newtheorem{theorem}{Theorem}[section]
\newtheorem{lemma}[theorem]{Lemma}
\newtheorem{corollary}[theorem]{Corollary}
\newtheorem{proposition}[theorem]{Proposition}

\theoremstyle{definition}

\newtheorem{conjecture}[theorem]{Conjecture}
\newtheorem{remark}[theorem]{Remark}

\newtheorem{example}[theorem]{Example}
\newtheorem{exercise}[theorem]{Exercise}

\definecolor{myyellow}{RGB}{255,241,172}
\definecolor{myteal}{RGB}{250,227,217}
\definecolor{mygreen}{RGB}{187,222,214}
\definecolor{myblue}{RGB}{138,198,209}

\newcommand{\p}{\mathfrak{p}}
\newcommand{\Frob}{\text{Frob}}
\newcommand{\Gal}{\text{Gal}}
\newcommand{\ind}{\text{ind}}

\begin{document}
	
	\title{Moments of ideal class counting functions}
	
	\author[Kam Cheong Au]{Kam Cheong Au}
	
	\address{Rheinische Friedrich-Wilhelms-Universität Bonn \\ Mathematical Institute \\ 53115 Bonn, Germany} 
	
	\email{s6kmauuu@uni-bonn.de}
	

	\begin{abstract} We consider the counting function of ideals in a given ideal class of a number field of degree $d$. This describes, at least conjecturally, the Fourier coefficients of an automorphic form on $\text{GL}(d)$, typically not a Hecke eigenform and not cuspidal. We compute its moments, and also investigate the moments of the corresponding cuspidal projection.

		%Motivated by the moments of positive-definite quadratic form counting function, we look at the moments of ideal class counting function on general number fields. Moreover, we define a closely and naturally associated "cuspdial" version of it, and investigate its moments.  
	\end{abstract}
	
	\maketitle
	
	\section{Introduction}
	Given a positive definite binary integral quadratic form $g(x,y)$ of discriminant $D$, the properties of the numbers $r_g(n) = \#\{(x,y)\in \mathbb{Z}^2 \mid g(x,y)=n\}$ have been intensively investigated. The Dirichlet series $\sum_{n\geq 1} r_g(n) n^{-s}$ is the $L$-function of an elliptic modular form, the moments $$\sum_{n\leq x} r_g(n)^{2\beta} \asymp x (\log x)^{2^{2\beta-1}-1}, \qquad \beta\in \mathbb{R}^{>0},$$
	have been quite well-studied (for example \cite{toma2020estimates}, \cite{rankin1983sums}, \cite{rankin1985sums} and \cite{odoni1991problem}). \par One can define, in a quite natural way, a closely related quantity $r_{\text{cusp},g}(n)$ such that $\sum_{n\geq 1}r_{\text{cusp},g}(n)n^{-s}$ is now the $L$-function of a cuspidal modular form. 
	One expects a slower growth of the following moment
	$$\sum_{n\leq x} |r_{\text{cusp},g}(n)|^{2\beta} \asymp x (\log x)^A.$$
	It is important to note that neither $r_g(n)$ nor $r_{\text{cusp},g}(n)$ are in general Fourier coefficients of a Hecke eigenform. Therefore the exponent $A$ cannot be predicted from a suitable Sato-Tate law and is an interesting quantity to consider. Blomer \cite{blomer2004cusp} gives a formula for $A$ in when $D$ is fundamental; we will derive a similar formula for non-fundamental $D$ as an application of our general framework to be described below. \\[0.02in]
	
	Because $r_g(n)$ is essentially number of integral ideals of norm $n$ in a given ideal class, one can extend this investigation to general number fields. Let $K$ be a number field of degree $d$, $\mathfrak{A}$ an ideal class of $K$, $a(\mathfrak{A},n)$ the number of integral ideals in the class $\mathfrak{A}$ with norm $n$. We first consider the moment
	$$\sum_{n\leq x} a(\mathfrak{A},n)^{2\beta} \asymp x (\log x)^{A_1}.$$
	We will give a formula for exponent $A_1$ (Theorem \ref{full_ideal_class_moment}). When $K/\mathbb{Q}$ is Galois, it turns out $A_1 = d^{2\beta-1}-1$; the non-Galois case is more complicated. \par
	
	Being a linear combination of $L$-functions of automorphic representations of $\text{GL}(1)/K$, one could perform (at least conjecturally) automorphic induction to produce representations of $\text{GL}(d)/\mathbb{Q}$, where it makes sense to talk about the cusp space. We study this phenomenon in the second section. We first give a natural criterion on which the corresponding Artin $L$-function should come from a cuspidal automorphic representation (Proposition \ref{cusp_criterion_prop}). This criterion agrees with proven cases of automorphic induction \cite[Section~3.6]{arthur1989simple}. \par
	
	In general number fields, we can again then define a quantity $a_\text{cusp}(\mathfrak{A},n)$ that generalizes $r_{\text{cusp},g}(n)$, and $\sum_{n\geq 1} a_\text{cusp}(\mathfrak{A},n)n^{-s}$ is a linear combination of $L$-functions of cuspidal automorphic forms. Consider the same moment problem:
	$$\sum_{n\leq x} |a_\text{cusp}(\mathfrak{A},n)|^{2\beta} \asymp x (\log x)^{A_2}.$$
	We will prove (in majority of cases: Proposition \ref{cusporderintegral}, \ref{cusporderreal}) a formula for $A_2$. It will involve finer properties of the number field $K$. \par
	
	In our investigation, the naive language of quadratic form no longer suffices. Instead, we employ extensively the language of (finite group) representations. Not only it allows us to state certain formulas more concisely, but also provides more transparent proofs for the known quadratic cases. \\[0.02in]
	
	We work mainly with Dirichlet series, the asymptotic of the partial sum is then obtained via Tauberian theorems \cite[Chapter~5.3]{tenenbaum2015introduction}. 
	\\[0.02in]
	
	Most of the works below are extracted from the author's Master thesis (at University of Bonn) under the same title. 
	~\\[0.05in]

	\section{Full partial zeta series}
	
	\subsection{General considerations}
	In this section, we outline the framework we will be working throughout this article and then prove several important results. We will temporarily not mention objects related to ideal classes. Doing so will allow us to concentrate on the more important representation-theoretic details. \par
	
	We fix the following notations which will be used in this section. Let $L/K$ be an abelian extension of number fields, $M$ be an extension of $L$ such that $M/\mathbb{Q}$ is Galois ($M$ can be taken, for example, to be the Galois closure of $L$ over $\mathbb{Q}$). Denote $$G = \text{Gal}(M/\mathbb{Q})\quad H = \text{Gal}(M/K) \quad N = \text{Gal}(L/K)
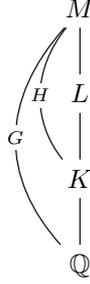
$$
	
	$$\begin{tikzcd}
		M  \\
		L \\
		K \\
		{\mathbb{Q}} 
		\arrow[no head, from=1-1, to=2-1]
		\arrow[no head, from=2-1, to=3-1]
		\arrow[no head, from=3-1, to=4-1]
		\arrow["H"{description}, curve={height=18pt}, no head, from=1-1, to=3-1]
		\arrow["G"{description}, curve={height=30pt}, no head, from=1-1, to=4-1]
	\end{tikzcd}$$
	
	Let $\chi: N\to \mathbb{C}^\times$ be a (one-dimensional) character of $N$, denote $\widetilde{\chi} : H \to \mathbb{C}^\times$ to be the one-dimensional character of $H$ by composing $\chi$ with the natural map $H\to N$ with $\chi$, let $\widetilde{\chi}^\ind$ be the induced character of $\widetilde{\chi}$ from $H$ to $G$ (it might not be multiplicative). We will sometimes abuse notations by using same symbols for characters and their associated representations. For two characters $\rho_1,\rho_2$ of $G$, write
	$$\langle \rho_1,\rho_2\rangle_G = \frac{1}{|G|} \sum_{g\in G} \rho_1(g) \overline{\rho_2(g)}$$
	as the usual inner product between characters. 
	
	Denote the Artin $L$-function associated with a representation $\rho$ of $G$ as $L(\rho,s)$, also let $a(\rho,n)$ be defined by
	$$L(\rho,s) = \sum_{n\geq 1} \frac{a(\rho,n)}{n^s}.$$
	
	\begin{lemma}
		We have $L(\chi,s) = L(\widetilde{\chi},s) = L(\widetilde{\chi}^\ind,s)$, that is 
		$$a(\chi,n) = a(\widetilde{\chi},n) = a(\widetilde{\chi}^\ind,n).$$
	\end{lemma}
	\begin{proof}
		This follows immediately from the fact that Artin $L$-function is invariant under restriction and induction. (See \cite[Chapter~7]{neukirch2006algebraische})
	\end{proof}
	
	For a prime $p$ in $\mathbb{Q}$, unramified in $M$, recall that the Frobenius $\Frob(p)$ is a well-defined conjugacy class in $G$: take any $\mathfrak{P}$ in $M$ lying over $p$, then $\Frob(p)$ is the conjugacy class containing $\Frob(\mathfrak{P}/p)$. 
	
	\begin{lemma}
		For a prime $p$ unramified in $\mathbb{Q}$, let $C$ be the corresponding conjugacy class of the Frobenius automorphism $\Frob(p)$ in $G$, then
		$$a(\widetilde{\chi}^\ind,p) = \widetilde{\chi}^\ind (C).$$
	\end{lemma}
	\begin{proof}
		This follows immediately by comparing the first order expansion of determinant: $$\frac{1}{\det(1-\widetilde{\chi}^\ind(\Frob p) p^{-s})} = 1 + \frac{\widetilde{\chi}^\ind(\Frob p)}{p^s} + \cdots$$
	\end{proof}
	
	For real $\beta>0$, we are interested in the behaviour of the series $\sum_{n\leq x} |a(\chi,n)|^{2\beta} n^{-s} $ near $s=1$. 
	
	\begin{proposition}\label{characterorder}
		Let $\chi$ be a $1$-dimensional character of $N$, then as $s\to 1$ from the right, we have, for some $c\neq 0$, 
		$$\sum_{n\geq 1} \frac{|a(\chi,n)|^{2\beta}}{n^s} \sim c (s-1)^{-\varrho(\chi,\beta)}$$
		where $$\varrho(\chi,\beta) = \frac{1}{|G|} \sum_{g\in G} |\widetilde{\chi}^\ind(g)|^{2\beta}$$
	\end{proposition}
	\begin{proof}
		We write $A\approx B$ if $A/B$ is holomorphic on $\Re(s)>1/2$ without zero. Then we have
		$$\sum_{n=1}^\infty \frac{|a(\chi,n)|^{2\beta}}{n^s} = \prod_p \left(1+\frac{|a(\chi,p)|^{2\beta}}{p^s} +\frac{|a(\chi,p^2)|^{2\beta}}{p^{2s}}+\cdots \right) \approx \prod_p \left(1+\frac{|a(\chi,p)|^{2\beta}}{p^s} \right).$$
		Here last step is valid since $a(\chi,n) = O(n^\varepsilon)$ for any $\varepsilon>0$. We can neglect those finitely many ramified primes. Let $C$ denote a conjugacy class in $G$, above is
		$$\prod_{C\subset G} \prod_{\Frob{p} = C} \left(1+ \frac{|\widetilde{\chi}^\ind(C)|^{2\beta}}{p^s}\right) \approx \prod_{C\subset G} \prod_{\Frob{p} = C} \left(1+ \frac{1}{p^s}\right)^{|\widetilde{\chi}^\ind(C)|^{2\beta}}.$$
		From Chebotarev's density theorem, we know $\prod_{\Frob{p} = C} (1+ \frac{1}{p^s})$ has polar density $|C|/|G|$. So we have $$\varrho(\chi,\beta) = \sum_{C\subset G} \frac{|C|}{|G|} |\widetilde{\chi}^\ind(C)|^{2\beta} = \frac{1}{|G|} \sum_{g\in G} |\widetilde{\chi}^\ind(g)|^{2\beta}.$$
	\end{proof}
	
	\begin{corollary}
		For $\beta>0$, there exists a constant $c\neq 0$ such that
		$$\sum_{n\leq x} |a(\chi,n)|^{2\beta} \sim c x (\log x)^{\varrho(\chi,\beta)-1}.$$
	\end{corollary}
	\begin{proof}
		This follows from the above after using a version of Tauberian theorems. See \cite[Theorem~4]{delange1954generalisation} or \cite[Chapter~5.3]{tenenbaum2015introduction}.
	\end{proof}
	
	We record an observation that will be used later:
	\begin{lemma}\label{meromorphic_cont}
		Let $\rho_1,\cdots,\rho_k$ be representations of $G$, then $$\sum_{n\geq 1}\frac{a(\rho_1,n)\cdots a(\rho_k,n)}{n^s},$$
		originally convergent for $\Re(s)>1$, has meromorphic continuation to $\Re(s) > 1/2$.
	\end{lemma}
	\begin{proof}
		Write $A\approx B$ if $A/B$ is holomorphic on $\Re(s)>1/2$ without zeroes. We can neglect those finitely ramified primes. Then 
		$$\sum_{n\geq 1}\frac{a(\rho_1,n)\cdots a(\rho_k,n)}{n^s} = \prod_p \left(1+\frac{a(\rho_1,p)\cdots a(\rho_k,p)}{p^s} + \frac{a(\rho_1,p^2)\cdots a(\rho_k,p^2)}{p^{2s}} + \cdots \right) \approx \prod_p \left(1+\frac{a(\rho_1,p)\cdots a(\rho_k,p)}{p^s} 
		\right).$$
		For unramified $p$, one has $a(\rho_1,p)a(\rho_2,p) = a(\rho_1\otimes \rho_2,p)$ with $\otimes$ meaning tensor product representation. Therefore above equals $$\prod_p \left(1+ \frac{a(\rho_1\otimes \cdots \otimes \rho_k,p)}{p^s}\right) \approx \prod_p \left(1+ \frac{a(\rho_1\otimes \cdots \otimes \rho_k,p)}{p^s} + \frac{a(\rho_1\otimes \cdots \otimes \rho_k,p^2)}{p^{2s}}+\cdots \right) \approx L(\rho_1\otimes \cdots \otimes \rho_k,s),$$
		this is the Artin $L$-function of another representation of $G$, which is known to admit meromorphic continuation.
	\end{proof}
	
	Using the main conclusion of  \cite{draxl1971funktionen}, it should not be too hard to prove that the Dirichlet series in the above lemma has meromorphic continuation to all of $\mathbb{C}$, but we will not need this fact. ~\\[0.02in]
	
	Denote $$\varrho(\beta) = \max_{\chi \in \widehat{N}} \varrho(\chi,\beta).$$
	Here the maximum is taken over all (1-dimensional) character of $N$. Let $\mathfrak{X} \subset \widehat{N}$ such that the maximum is attained, we will see in proof of following lemma that $\mathfrak{X}$ is independent of $\beta$. 
	
	\begin{lemma}\label{extremealcharlemma}
		(i) We have $$\varrho(\beta) = \frac{1}{|G|} \sum_{C\subset G} |C| \left( \frac{[G:H] |H\cap C|}{|C|}\right)^{2\beta},$$
		here we are summing over all conjugacy classes of $G$.\\
		(ii) Let $\chi \in \mathfrak{X}$. Then for each $C \subset G$, $\widetilde{\chi}$ is constant on $C\cap H$, and $\widetilde{\chi}^\ind (g)$ is of the form $\mathbb{Z}^{\geq 0}$ times a root of unity.\\
		(iii) Given $\sigma\in N$, there exists a conjugacy class $C\subset G$ such that $\overline{\chi(\sigma)} \widetilde{\chi}^\ind (C) \geq 0$ for all $\chi \in \mathfrak{X}$.
	\end{lemma}
	\begin{proof}
		Let $C$ be a conjugacy class of $G$, $\phi$ be the indicator function on $C$. We have $$\langle \widetilde{\chi}^\ind,\phi\rangle_G = \frac{1}{|G|}\sum_{g\in C} \widetilde{\chi}^\ind(g).$$
		On the other hand, by Frobenius reciporcity,
		$$\langle \widetilde{\chi}^\ind,\phi\rangle_G = \langle \widetilde{\chi},\text{Res }\phi\rangle_H = \frac{1}{|H|} \sum_{h\in H\cap C} \widetilde{\chi}(h),$$
		thus $$\widetilde{\chi}^\ind(C) = \frac{[G:H]}{|C|} \sum_{h\in H\cap C}\widetilde{\chi}(h).$$ Hence $$\varrho(\chi,\beta) = \frac{1}{|G|} \sum_{C\subset G} \frac{|C|}{|G|} |\widetilde{\chi}^\ind(C)|^{2\beta}$$ is maximized when $\chi$ is the trivial character on $N$. For any other character $\chi$ also maximizing it, we must have $|\widetilde{\chi}^\ind(C)| = |1^\ind(C)|$ for all $C$, implying tha $\widetilde{\chi}(h)$ is constant on $h\in H\cap C$, and the fact that $\widetilde{\chi}(h)$ is a root of unity implies $\widetilde{\chi}^\ind(C)$ is a positive integer times roots of unity. \par
		For the last assertion, choose any lift $\sigma' \in H$ of $\sigma\in N$, let $C$ be the conjugacy class of $G$ containing $\sigma'$, then this $C$ satisfies the condition. 
	\end{proof}
	
	Consider the following partial zeta coefficient: for $\sigma\in N$, let\footnote{we used the notation $a(\cdot,n)$ in two different ways: if $\sigma \in N$, $a(\sigma,n)$ is defined via the next displayed equation; if $\chi \in \widehat{N}$, $a(\chi,n)$ is the coefficient of the Artin $L$-function $L(\chi,s)$. The symbol $\sigma$ will always denote an element in $N$; $\chi$ will always denote an element in $\widehat{N}$.}$$a(\sigma,n) = \frac{1}{|N|} \sum_{\chi \in \widehat{N}} \bar{\chi}(\sigma) a(\chi,n).$$

	Let $S$ be a finite collection of finite places in $\mathbb{Q}$, i.e. a finite set of prime numbers, write $(n,S) = 1$ if $n$ is not divisible by each element in $S$. 
	
	\begin{lemma}\label{partial_zeta_interpretation_lemma}
		For any $\sigma\in N$ and $n \in \mathbb{Z}^{\geq 1}$, one has $a(\sigma,n) \in \mathbb{Q}$. Let $S$ the be finite primes in $\mathbb{Q}$ which lies below primes that are ramifieid in $L/K$, then we also have $a(\sigma,n) \in \mathbb{Z}^{\geq 0}$ for $(n,S)=1$.
	\end{lemma}
	\begin{proof}
		Let $\tau \in \Gal(\overline{\mathbb{Q}}/\mathbb{Q})$, denote $\chi^\tau$ denote the character by composing the value of $\chi$ with $\tau$. Then $$a(\sigma,n)^\tau = \frac{1}{|N|}\sum_{\chi \in \widehat{N}}\overline{\chi^\tau}(\sigma) a(\chi,n)^\tau = \frac{1}{|N|}\sum_{\chi \in \widehat{N}}\overline{\chi^\tau}(\sigma) a(\chi^\tau,n) = \frac{1}{|N|}\sum_{\chi \in \widehat{N}}\overline{\chi}(\sigma) a(\chi,n) = a(\sigma,n),$$
		this holds for all $\tau \in \Gal(\overline{\mathbb{Q}}/\mathbb{Q})$, so $a(\sigma,n)\in \mathbb{Q}$. Let $S'$ be primes in $K$ lying over $S \subset \mathbb{Q}$. Note that $$\sum_{n\geq 1} \frac{a(\chi,n)}{n^s} = \prod_{\p \in S'} L_\p(\chi,s) \prod_{\p\notin S'} \frac{1}{1-\chi(\p)(N\p)^{-s}},$$
		here the Euler factors $L_\p(\chi,s)$ for ramified primes lack a uniform description; but for $(n,S)=1$, it does not affect value of $a(\sigma,n)$, so $a(\sigma,n)$ is non-negative integer from character orthogonality.
	\end{proof}

	Our first major result is
	\begin{theorem}\label{partialorder}
		Let $\beta>0, \sigma\in N$, then as $s\to 1$ from the right, we have
		$$0< \liminf_{s\to 1} (s-1)^{\varrho(\beta)} \sum_{n\geq 1} \frac{|a(\sigma,n)|^{2\beta}}{n^s} \leq \limsup_{s\to 1} (s-1)^{\varrho(\beta)}\sum_{n\geq 1} \frac{|a(\sigma,n)|^{2\beta}}{n^s} < \infty.$$
	\end{theorem}
	\begin{proof}[Proof of Theorem \ref{partialorder}]
		The fact that $\limsup < \infty$ is evident since $$|a(\sigma,n)|^{2\beta} \ll \sum_\chi |a(\chi,n)|^{2\beta}$$
		and Dirichlet series of each term has exponent $\leq \varrho(\beta)$. \par
		Proving $\liminf \neq 0$ is more involved. Abbreviate $\widetilde{\chi}^\ind(C) := a_{C,\chi}$ where $C$ is a conjugacy class of $G$. For $\chi \in \mathfrak{X}$, by the above lemma, we can choose positive integer $l$ such that $a_{C,\chi}^l \geq 0$ for all $C$. Write $\mu$ as a primitive $l$-th root of unity. For some $\rho_i(C)\in \mathbb{C}, i=0,\cdots,l-1$ to be fixed later, define
		$$f(\chi,C) = \sum_{i=0}^{l-1} \rho_i(C) \prod_{\Frob p \in C} \left(1+ \frac{\mu^i a_{C,\chi}}{p^s} \right)$$ 
		also define $a'(\chi,n)$ and $a'(\sigma,n)$ via
		$$\sum_{n\geq 1} \frac{a'(\chi,n)}{n^s} := \prod_{C\subset G} f(\chi,C),$$
		$$a'(\sigma,n) = \frac{1}{|N|} \sum_{\chi\in \widehat{N}} \bar{\chi}(\sigma) a'(\chi,n).$$
		
		Note that the same prime $p$ never occurs in two $f(\chi,C_1)$ and $f(\chi,C_2)$ with $C_1\neq C_2$; as usual, we only focus on primes which are not ramified in $M$.
		\\[0.2in]
		We claim that we can choose $\rho_i(C) \neq 0$ independent of $\chi$ such that 
		\begin{enumerate}
			\item For all $j$ and $C$, $\sum_{i=0}^{l-1} \rho_i(C) \mu^{ij}$ equals either $0$ or $1$.
			\item $\bar{\chi}(\sigma)a'(\chi,n) \geq 0$ for all $n$ and $\chi \in \mathfrak{X}$. 
		\end{enumerate}
		
		Assuming above claim, let us prove $\liminf > 0$. (1) implies $a'(\chi,n) = a(\chi,n)$ for all $\chi$ or $a'(\chi,n) = 0$ for all $\chi$, therefore $|a(\sigma,n)| \geq |a'(\sigma,n)|$, so it suffices to prove the assertion for $\sum_{n\geq 1} |a'(\sigma,n)|^{2\beta}n^{-s}$. Since
		$$\begin{aligned}|a'(\sigma,n)|^{2\beta} &\gg \left|\sum_{\chi \in \mathfrak{X}} \bar{\chi}(\sigma)a'(\chi,n) \right|^{2\beta} - \sum_{\chi \notin \mathfrak{X}} |a'(\chi,n)|^{2\beta} \\
			&\gg \sum_{\chi \in \mathfrak{X}} |a'(\chi,n)|^{2\beta} - \sum_{\chi \notin \mathfrak{X}} |a'(\chi,n)|^{2\beta} \qquad \text{by condition (2)},\end{aligned}$$
		condition (1) also implies 
		\begin{equation}\tag{*}\sum_{n\geq 1} \frac{|a'(\chi,n)|^{2\beta}}{n^s} = \prod_{C\subset G} \left( \sum_{i=0}^{l-1} \rho_i(C) \prod_{\Frob p \in C} (1+ \frac{\mu^i |a_{C,\chi}|^{2\beta}}{p^s}) \right). \end{equation}
		To see why $(*)$ is true, we compute, for example coefficient of $p_1^{-s}p_2^{-s}$ for $p_1\neq p_2$, if $\Frob p_1 = \Frob p_2 = C$, then $$a'(\chi,p_1p_2) = \left(\sum_{i=0}^{l-1} \rho_i(C) \mu^{2i} \right)a_{C,\chi}^2.$$
		if $\Frob p_1 = C_1 \neq \Frob p_2 = C_2$, then
		$$a'(\chi,p_1p_2) = \left(\sum_{i=0}^{l-1} \rho_i(C_1) \mu^{i} \right)\left(\sum_{i=0}^{l-1} \rho_i(C_2) \mu^{i}\right) a_{C_1,\chi} a_{C_2,\chi}.$$
		since the numbers in parenthesis are either $0$ or $1$, we can distribute $|\cdot|^{2\beta}$ on both sides while keeping equalities. Above formulas generalize when $n=p_1\cdots p_k$, so $(*)$ indeed holds. \par
		
		In $(*)$, the $$\prod_{\Frob p \in C} (1+ \frac{\mu^i |a_{C,\chi}|^{2\beta}}{p^s}) \asymp (s-1)^{-\mu^i |a_{C,\chi}|^{2\beta} |C|/|G|},\qquad s\to 1.$$ As $i$ varies from $0$ and $l-1$, the dominant term come from those such that $\mu^i$ has the largest real part, i.e. $i=0$, so $$\sum_{i=0}^{l-1} \rho_i(C) \prod_{\Frob p \in C} \left(1+ \frac{\mu^i |a_{C,\chi}|^{2\beta}}{p^s} \right) \asymp (s-1)^{-|a_{C,\chi}|^{2\beta} |C|/|G|}$$ (here we also used $\rho_i(C)\neq 0$). Multiplying over all $C\subset G$ gives $\sum_{n\geq 1} |a'(\chi,n)|^{2\beta} n^{-s}$ is $\asymp (s-1)^{-\varrho(\chi,\beta)}$, which is $(s-1)^{-\varrho(\beta)}$ when $\chi \in \mathfrak{X}$. 
		Therefore $$\sum_{n\geq 1} \frac{|a'(\sigma,n)|^{2\beta}}{n^s} \gg (s-1)^{-\varrho(\beta)} - \sum_{\chi \notin \mathfrak{X}} \sum_{n\geq 1} \frac{|a'(\chi,n)|^{2\beta}}{n^s}.$$
		by definition of $\mathfrak{X}$, all terms on the right has pole order $<\varrho(\beta)$, so above is still $\gg (s-1)^{-\varrho(\beta)}$, proving $\liminf$ is positive, assuming two conditions above. 
		\\[0.2in]
		Now we explain how to choose $\rho_i(C)$ achieving these criteria. We make
		\begin{equation}\label{rhochoice1}\rho_i(C) = l^{-1} \quad \forall i \qquad   \qquad\text{ or } \qquad \rho_i(C) = l^{-1} \mu^{-1-i} \quad \forall i.\end{equation}
		Obviously $(1)$ is satisfied. For $(2)$, the lemma above says there exists conjugacy class $C_0 \subset G$ such that $\bar{\chi}(\sigma) a_{C_0,\chi} = \bar{\chi}(\sigma) \widetilde{\chi}^\ind(C_0) \geq 0$ for all $\chi \in \mathfrak{X}$. For $C\neq C_0$, we pick the first possibility in \ref{rhochoice1}, giving
		$$f(\chi,C) = 1 + \sum_{p_i\in \Frob C, p_i\neq p_j} \frac{a_{C_0,\chi}^l}{p_1^s \cdots p_l^s} + \sum_{p_i\in \Frob C, p_i\neq p_j} \frac{a_{C,\chi}^{2l}}{p_1^s \cdots p_{2l}^s} + \cdots,  \qquad C\neq C_0,$$
		which has non-negative coefficients (by our choice of $l$); for $C = C_0$, we pick the second possibility in (\ref{rhochoice1}), giving
		$$f(\chi,C_0) = \sum_{p_1\in \Frob C_0} \frac{a_{C,\chi}}{p_1^s} + \sum_{p_i\in \Frob C_0, p_i\neq p_j} \frac{a_{C_0,\chi}^{l+1}}{p_1^s \cdots p_{l+1}^s} + \sum_{p_i\in \Frob C_0, p_i\neq p_j} \frac{a_{C_0,\chi}^{2l+1}}{p_1^s \cdots p_{2l+1}^s} + \cdots,$$
		so $\bar{\chi}(\sigma)f(\chi,C_0)$ has non-negative coefficients. Then
		$$\sum_{n\geq 1} \frac{\bar{\chi}(\sigma)a'(\chi,n)}{n^s} = \bar{\chi}(\sigma) f(\chi,C_0) \prod_{C\neq C_0} f(\chi,C)$$ also has non-negative coefficients, which is $(2)$. 
	\end{proof}
	
	The above proof is an extended version of that used by Blomer in \cite{blomer2004cusp}, where it is essentially our case specialized to $l=2$.
	
	\begin{remark}Examining the above proof shows that we do not need to sum over all $n\geq 1$ to make $\liminf >0$. In fact, for any finite set of primes $S$ in $\mathbb{Q}$, replacing $\sum_{n\geq 1}|a(\sigma,n)|^{2\beta} n^{-s}$ with $$\sum_{\substack{(n,S)=1 \\ n \text{ squarefree}}} \frac{|a(\sigma,n)|^{2\beta}}{n^s}$$
		produces exactly the same qualitative behaviour, and the statement of above theorem remains unchanged. We only need to modify the proof to exclude primes $p \in S$ in the expression $f(\chi,C)$ above.
	\end{remark}
	
	For $\chi_1,\cdots,\chi_k \in \widehat{N}$, imitating the proof of Proposition \ref{characterorder}, one shows easily the pole of $$\sum_{n\geq 1} \frac{a(\chi_1,n)\cdots a(\chi_k,n)}{n^s}$$ at $s=1$ has order $$\varrho(\chi_1,\cdots,\chi_n) := \frac{1}{|G|} \sum_{g\in G} \widetilde{\chi_1}^\ind(g) \cdots \widetilde{\chi_k}^\ind(g) \in  \mathbb{Z}^{\geq 0}.$$
	
	\begin{lemma}\label{noncusp_max_lemma}
		Let $\chi_1,\cdots,\chi_k \in \widehat{N}$. If $\varrho(\chi_1,\cdots,\chi_k) = \varrho(k/2)$, then each $\chi_i \in \mathfrak{X}$ and $\chi_1\cdots \chi_k = 1$. 
	\end{lemma}
	\begin{proof}
		Hölder's inequality implies
		$$\varrho(\chi_1,\cdots,\chi_k) \leq \left(\frac{1}{|G|} \sum_{g\in G} |\widetilde{\chi_1}^\ind(g)|^k \right)^{1/k}\cdots \left(\frac{1}{|G|} \sum_{g\in G} |\widetilde{\chi_k}^\ind(g)|^k \right)^{1/k} = \varrho(\chi_1,k/2)^{1/k} \cdots \varrho(\chi_k,k/2)^{1/k}$$
		which is $\leq \varrho(k/2)$. Equality holds if and only if each term is equals $\varrho(k/2)$, that is $\chi_i \in \mathfrak{X}$. 
		For the second assertion, from Lemma \ref{extremealcharlemma}, we know that $\chi_i$ is constant on each $C\cap H$ if it's non-empty (otherwise 0), and we have $$\widetilde{\chi_i}^\ind(C) = \frac{[G:H]|C\cap H|}{|C|} \widetilde{\chi_i}(H\cap C),$$
		hence $$\varrho(\chi_1,\cdots,\chi_k) = \frac{1}{|G|}\sum_{C\subset G, C\cap H\neq \varnothing} \left( \frac{[G:H]|C\cap H|}{|C|}\right)^k \widetilde{\chi_1}(C\cap H) \cdots \widetilde{\chi_k}(C\cap H),$$
		so if this equals $\varrho(k/2)$, $\widetilde{\chi_1}(C\cap H) \cdots \widetilde{\chi_k}(C\cap H)$ must be $1$ for all conjugacy classes $C$ of $G$ that intersect $H$, so $\chi_1\cdots \chi_k = 1$. 
	\end{proof}
	
	\begin{corollary}\label{ideal_counting_moment_corollary}
		Let $S$ be finite primes in $\mathbb{Q}$ which lie below primes that are ramifieid in $L/K$. Then for any integer $k\geq 1$, there exists constant $C_{k,S} > 0$, independent of $\sigma$, such that
		$$\sum_{n\geq 1, (n,S)=1} \frac{a(\sigma,n)^k}{n^s} \sim C_{k,S} (s-1)^{\varrho(k/2)}, \qquad s\to 1,$$
		and
		$$\sum_{n\leq x, (n,S)=1} a(\sigma,n)^k \sim \frac{C_{k,S}}{\varrho(k/2)!} x (\log x)^{\varrho(k/2)-1}.$$
	\end{corollary}
	\begin{proof}
		Since $(n,S)=1$, each $a(\sigma,n)$ is non-negative, so we can remove the absolute value in Theorem \ref{partialorder}. On the other hand, from Lemma \ref{meromorphic_cont}, we know that $f(s)=\sum_{n\geq 1, (n,S)=1} a(\sigma,n)^k n^{-s}$ has meromorphic extension to $\Re(s)>1/2$, then above theorem implies that it has pole of order exactly $\varrho(k/2)$ there. \par 
		The constant $C_{k,S}$ is the leading coefficient of $f(s)$, it remains to prove it is independent of $\sigma$. We have
		$$f(s) = \frac{1}{|N|^k} \sum_{\chi_1,\cdots,\chi_k \in \widehat{N}} \overline{(\chi_1 \cdots \chi_k)(\sigma)} \sum_{n\geq 1,(n,S)=1}a(\chi_1,n)\cdots a(\chi_k,n) n^{-s}.$$
		Here $\sigma$ only appears as $(\chi_1 \cdots \chi_k)(\sigma)$. The inner Dirichlet series has pole of order $\varrho(\chi_1,\cdots,\chi_n)$, if it contributes to leading coefficient, the above lemma forces $\chi_1\cdots \chi_k = 1$, so $(\chi_1 \cdots \chi_k)(\sigma) = 1$, thus it is independent of $\sigma$.
	\end{proof}
	
	When $\beta\in \mathbb{Z}/2$, we see that $\liminf$ and $\limsup$ in Theorem \ref{partialorder} are equal, it is very natural to expect this also holds for $\beta>0$:
	\begin{conjecture}
		For real $\beta>0$, the $\liminf$ and $\limsup$ in Theorem \ref{partialorder} are equal, and they are independent of $\sigma\in N$. 
	\end{conjecture}
	
	\subsection{Ideal class counting function}
	
	In this section, we apply the results proved above to ideal classes in number field. 
	
	Let $K$ be a number field, $L$ its Hilbert class field, $M$ be the Galois closure of $L/\mathbb{Q}$. Using all previous notations, we see that $N$ is isomorphic to the ideal class group: $\Gal(L/K) = N$ corresponds to ideal classes of $K$ under inverse of Artin map. Moreover, $$\sum_{n\geq 1} \frac{a(\chi,n)}{n^s} = L(\chi,s) = \sum_{I \subset \mathcal{O}_K} \frac{\chi(I)}{N(I)^s}$$
	is the $L$-series of an ideal class group character and $$a(\sigma,n) = \frac{1}{|N|} \sum_{\chi \in \widehat{N}} \bar{\chi}(\sigma) a(\chi,n)$$ counts the number of integral ideals of norm $n$ in the class $\sigma$. 
	
	\begin{theorem}\label{full_ideal_class_moment}
		Let $K$ be a number field, let $\mathfrak{A}$ denote an ideal class of $K$, $a(\mathfrak{A},n)$ denote the number of integral ideals in this class with norm $n$. For $k\geq 1$ positive integer, there exists constants $C_k>0$, independent of $\mathfrak{A}$, such that as $s\to 1$,
		$$\sum_{n\geq 1} \frac{a(\mathfrak{A},n)^k}{n^s} \sim C_k (s-1)^{\varrho(k/2)} \qquad s\to 1,$$
		and
		$$\sum_{n\leq x} a(\mathfrak{A},n)^k \sim \frac{C_k}{\varrho(k/2)!} x (\log x)^{\varrho(k/2)-1},$$
		with the positive integer $\varrho(k/2)$ computed using notations in previous section. Moreover, for general real $\beta>0$, 
		$$x (\log x)^{\varrho(\beta)-1} \ll \sum_{n\leq x} a(\mathfrak{A},n)^{2\beta} \ll x (\log x)^{\varrho(\beta)-1}.$$
	\end{theorem}
	\begin{proof}
		This is a direct translation of results in previous section. Note that we can take $S = \varnothing$ since $L/K$ is unramified. 
	\end{proof}
	
	\begin{remark}
		When $k=1$, one easily computes $\varrho(k/2) = 1$, so $$\sum_{n\leq x} a(\mathfrak{A},n) \sim C_1 x.$$
		This is actually a special case of the following more general result (\cite[Chapter~6]{marcus1977number}):
		$$\sum_{n\leq x} a(\mathfrak{A},n) = \kappa_K x + O(x^{1-1/[K:\mathbb{Q}]}),$$
		where $\kappa_K$ is the residue of Dedekind zeta function of at $s=1$ divided by class number $h_K$ of $K$. 
	\end{remark}
	
	% \begin{remark}
		% When $k\geq 2$, the number $C_k$ is a "natural number" (i.e. it comes from values of $GL(n)$ $L$-functions at integer points) only if $K$ is quadratic and $k=2$. We will elucidate its nature in the last section. 
		% \end{remark}
	
	The quantities $\varrho(\chi, \beta), \varrho(\beta)$ are in general difficult to compute as $[M:\mathbb{Q}]$ could be very large compared to $[K:\mathbb{Q}]$. There is a simpler formula when $K/\mathbb{Q}$ is assumed to be Galois.
	\begin{proposition}
		Let $K/\mathbb{Q}$ be Galois, $\chi$ a character of its ideal class group $\mathcal{C}$, $$\varrho(\chi,\beta) = \frac{1}{h_K [K:\mathbb{Q}]} \sum_{\mathfrak{A}\in \mathcal{C}} |l(\chi,\mathfrak{A})|^{2\beta},$$
		where $l(\chi,\mathfrak{A}) = \sum_{\tau \in \Gal(K/\mathbb{Q})} \chi(\mathfrak{A}^\tau)$. Moreover, for $\beta>0$, $$\varrho(\beta) = \max_{\chi} \varrho(\chi,\beta) = [K:\mathbb{Q}]^{2\beta-1}.$$
	\end{proposition}
	\begin{proof}
		When $K/\mathbb{Q}$ is Galois, $L/\mathbb{Q}$ is also Galois\footnote{If $K/\mathbb{Q}$ is Galois and $L$ Hilbert class field (= maximal unramified abelian extension) of $K$, then $L/\mathbb{Q}$ is also Galois. This is a nice exercise in algebraic number theory, for completeness we quickly recall the proof: we need to show for each $\sigma \in \Gal(\overline{\mathbb{Q}}/\mathbb{Q})$, $\sigma(L) = L$; since $L/K$ is abelian unramified, so is $\sigma(L)/\sigma(K) = \sigma(L)/K$ (we used here the assumption $K/\mathbb{Q}$ is Galois), but $L$ is maximal under this property, so $\sigma(L)\subset L$, repeating the argument with $\sigma$ replaced by $\sigma^{-1}$ shows the other inclusion, so $\sigma(L)=L$.}, so we can take $M=L$ and $N=H$, thus there is no difference between $\chi$ and $\widetilde{\chi}$. Also $H$ is a normal subgroup of $G$, its induced character then takes the following special form:
		$$\chi^\ind(g) = \begin{cases} 0 \qquad & g\notin H \\ \sum_{s\in G/H} \chi(s^{-1}gs) \qquad & g\in H \end{cases}$$
		and $$\varrho(\chi,\beta) = \frac{1}{|G|}\sum_{g\in G} |\chi^\ind(g)|^{2\beta} = \frac{1}{h_K [K:\mathbb{Q}]}\sum_{h\in H} |\chi^\ind(h)|^{2\beta}.$$
		If $h\in H$ correspond to an ideal class $\mathfrak{A}$, then as $s$ varies over $G/H$, $s^{-1} h s$ are exactly ideal classes $\mathfrak{A}^\tau$ with $\tau \in \Gal(K/\mathbb{Q})$, this gives the term $l(\chi,\mathfrak{A})$.\par
		For the formula of $\varrho(\beta)$, we know the maximum is attained when $\chi$ is trivial, and $l(1,\mathfrak{A}) = [K:\mathbb{Q}]$, so
		$$\varrho(\beta) = \frac{1}{h_K [K:\mathbb{Q}]} \sum_{h\in H} [K:\mathbb{Q}]^{2 \beta} = [K:\mathbb{Q}]^{2 \beta - 1}.$$
	\end{proof}
	
	\begin{example}\label{quad_example} If $K/\mathbb{Q}$ is quadratic, $\text{Gal}(K/\mathbb{Q}) = \{1,\sigma\}$, let $m$ be order of $\chi$, $\mathfrak{A}_j$ be the set of ideal classes $\mathfrak{A}$ such that $\chi(\mathfrak{A}) = e^{2\pi i j/m}$. Then $|\mathfrak{A}_j| = h/m$ and $$l(\chi,\mathfrak{A}) = \chi(\mathfrak{A}) + \chi(\mathfrak{A}^\sigma) = \chi(\mathfrak{A}) + \chi(\mathfrak{A}^{-1}) = 2\cos(\frac{2\pi j}{m}),  \qquad \mathfrak{A}\in \mathfrak{A}_j.$$
		So $$\varrho(\chi,\beta) = \frac{1}{2h} \sum_{j=0}^{m-1} \sum_{\mathfrak{A}\in \mathfrak{A}_j} |l(\chi,\mathfrak{A})|^{2\beta} = \frac{1}{2h} \sum_{j=0}^{m-1} \frac{h}{m} \left|2\cos(\frac{2\pi j}{m})\right|^{2\beta} = \frac{1}{2m} \sum_{j=0}^{m-1} \left|2\cos(\frac{2\pi j}{m})\right|^{2\beta},$$
		which recovers a formula in Blomer \cite{blomer2004cusp}.
	\end{example}
	
	\begin{example}\label{cubicexample1}
		Let $K = \mathbb{Q}[x]/(x^3-21x-28)$, it's a cubic Galois field, let $\sigma$ be a generator of Galois group. $K$ has class group $\mathcal{C}\cong \mathbb{Z}/3\mathbb{Z}$, 2 splits in $K$, let $\p$ be a prime lying above $2$, it generates the class group. Let $\chi \in \widehat{\mathcal{C}}$ be an ideal class group character, such that $\chi(\p) = e^{2\pi i /3}$. \par 
		Since Galois group act transitively on prime ideals lying over $2$, $(2) = \p \p^\sigma \p^{\sigma^2}$, let $\p^\sigma \sim \p^i$, then $1+i+i^2 \equiv 0 \pmod{3}$, so $i\equiv 0 \pmod{3}$, thus $\p^\sigma \sim \p$ are in the same ideal class. Therefore 
		$$l(\chi,1) = 3, \qquad l(\chi,\p) = 3\chi(\p) = 3e^{2\pi i /3}, \qquad l(\chi,\p^2) = 3e^{-2\pi i /3}.$$
		So for every $\chi \in \widehat{\mathcal{C}}$, $$\varrho(\chi,\beta) = \frac{1}{9} (3^{2\beta} + |3e^{2\pi i /3}|^{2\beta} + |3e^{-2\pi i /3}|^{2\beta}) = 3^{2\beta -1}.$$
		in this case, every character attains the maximum $\varrho(\beta)$.
	\end{example}
	
	\begin{example}\label{cubicexample2}
		Let $K = \mathbb{Q}[x]/(x^3-x^2-54x+169)$, it's the unique cubic subfield of $\mathbb{Q}(\zeta_{163})$, let $\sigma$ be a generator of Galois group. $K$ has class group $\mathcal{C}\cong \mathbb{Z}/2\mathbb{Z} \times \mathbb{Z}/2\mathbb{Z}$, 5 splits in $K$, let $\p_1, \p_2$ be two primes lying above $5$, they generate the class group. \par 
		Let $A\in \text{GL}_2(\mathbb{Z}/2\mathbb{Z})\cong S_3$ be such that
		$$\left(\begin{smallmatrix} \p_1 \\ \p_2\end{smallmatrix} \right)^\sigma = A\left(\begin{smallmatrix} \p_1 \\ \p_2\end{smallmatrix} \right).$$
		Because $\sigma$ has order $3$, $A$ must also has order $3$, we can assume it is $A = \left(\begin{smallmatrix} 0 & 1 \\ 1 & 1\end{smallmatrix} \right)$, then $\p_1^\sigma \sim \p_2, \p_2^\sigma \sim \p_1 \p_2$. Therefore $\Gal(K/\mathbb{Q})$ acts transitively on non-principal ideal classes. \par
		Hence for $\mathfrak{A}\neq 1$ in $\mathcal{C}$, we have
		$$l(\chi,\mathfrak{A}) = \sum_{I\in \mathfrak{A}} \chi(I) - 1 = \begin{cases}3 \quad &\chi = 1 \\ -1 \quad &\chi \neq 1\end{cases}.$$
		Conclusion:
		$$\varrho(\chi,\beta) = \begin{cases}3^{2\beta-1} & \chi = 1 \\ \frac{1}{12}(3^{2\beta} + 1^{2\beta} + 1^{2\beta} + 1^{2\beta}) = \frac{1}{4}(3^{2\beta-1}+1) &\chi\neq 1\end{cases}.$$
		in this case, only the trivial character attains maximum $\varrho(\beta)$.
	\end{example}
	
	When $K/\mathbb{Q}$ is Galois, we have $\varrho(\beta) = [K:\mathbb{Q}]^{2\beta-1}$, this is not the case when $K/\mathbb{Q}$ is non-Galois. 
	\begin{example}
		Let $K/\mathbb{Q}$ be a non-Galois cubic field of class number $1$ (so we can take $L=K$). Let $a(n)$ be the number of integral ideals of norm $n$ in $K$, we claim $$\sum_{n\leq x} a(n)^{2\beta} \asymp x(\log x)^{\varrho(\beta)} \qquad \text{with }\varrho(\beta) = (1+3^{2\beta-1})/2.$$
		Indeed, $G = S_3$ has three conjugacy classes, and we can assume $H= \{(12),1\}$, $[G:H] = 3$. Using the formula $$\varrho(\beta) = \frac{1}{|G|} \sum_{C\subset G} |C| \left( \frac{[G:H] |H\cap C|}{|C|}\right)^{2\beta},$$
		we obtain
		$$\varrho(\beta) = \frac{1}{6}\left[2(\frac{3\times 0}{2})^{2\beta} + 3(\frac{3\times 1}{3})^{2\beta} + 1(\frac{3\times 1}{1})^{2\beta}\right],$$
		as claimed.
	\end{example}

	Next we present a computational example where $K/\mathbb{Q}$ is non-Galois to illustrate in general finding $\varrho(\beta)$ is not trivial at all.

	\begin{example}
		Let $K = \mathbb{Q}[X]/(X^4+5X^2-X+1)$, its class number is $3$, the Galois group of $K$ is $S_4$. We find the formula for $\varrho(\beta)$ with help of PARI, MAGMA and GAP.\footnote{it seems MAGMA alone suffices for all computations below.}

		Hilbert class field $L$ of $K$ has degree 12 over $\mathbb{Q}$, we can compute an absolute defining polynomial of $L$ in PARI using \begin{verbatim}K = bnfinit(y^4+5*y^2-y+1); polredabs(bnrclassfield(K,,2))\end{verbatim} 
		which outputs $f(x):= x^{12}-x^{11}+x^{10}-x^9-x^7+x^6-4 x^5-3 x^4+3 x^3+7 x^2+4 x+1$. Using \begin{verbatim}
			L = nfinit(f); nfsubfields(L,4)
		\end{verbatim}
		which is used to find degree $4$ subfields of $L$, including $K$, it gives a polynomial $g \in \mathbb{Z}[x]$ such that for $y = g(x)$ and $L = \mathbb{Q}[x]/(f(x))$, we have $K = \mathbb{Q}(y)$. Explicitly, $g(x) = -3 x^{11}+4 x^{10}-4 x^9+4 x^8-x^7+3 x^6-4 x^5+13 x^4+5 x^3-12 x^2-17 x-4$. \\
		
		Next we compute the Galois closure $M$ of $L/\mathbb{Q}$, this is a very large field. We switch to MAGMA since PARI is not capable finding Galois information $f(x)$ due to its large degree. Inputting 
		\begin{verbatim}
			P< x >:=PolynomialAlgebra(Rationals());
			f:=x^12 - x^11 + x^10 - x^9 - x^7 + x^6 - 4*x^5 - 3*x^4 + 3*x^3 + 7*x^2 + 4*x + 1;
			G, R, :=GaloisGroup(f: Prime:= 1913);
			print(G);
		\end{verbatim}
		we will explain the choice \texttt{Prime:= 1913} later. \texttt{G} gives the Galois group of $f$ in terms of permutations of its 12 roots: $G = \Gal(M/\mathbb{Q}) \subset S_{12}$ is generated by permutations
		\begin{multline*}\{(2, 12)(4, 7)(8, 10), \quad
			(2, 3, 4)(5, 10, 7)(8, 11, 12), \quad
			(3, 5, 11)(4, 12, 10),\quad (2, 7, 8)(4, 12, 10)\\
			(1, 5)(2, 4)(3, 9)(6, 11)(7, 10)(8, 12), \quad
			(1, 7)(2, 9)(3, 4)(5, 10)(6, 8)(11, 12), \\
			(1, 6, 9)(2, 8, 7)(3, 11, 5)(4, 10, 12)\}.\end{multline*}
		which has order $648 = [M:\mathbb{Q}]$. 
		
		WLOG, we may identify $\Gal(M/L)$ as the stabilizer of point $1$. We still need to identify what is $H = \Gal(M/K)$ under the above representation of $G$. We can work in any field extension of $\mathbb{Q}$ which contains all roots of $f$, for example $\mathbb{C}$ or $\mathbb{Q}_p$, we choose to work with later. $\mathbb{Q}_p$ contains all roots of $f$ if $p$ is unramified and $f \equiv 0$ splits modulo $p$. Such $p$ has density $1/[M:\mathbb{Q}] = 1/648$ and $p=1913$ satisfies the condition. 
		
		Such an explicit correspondence can be retrieved using \texttt{print(R)} in MAGMA. It says the index $i$ corresponds to $x_i \in \mathbb{Q}_{1913}$;
		$$\begin{aligned}
			x_1 &= 181 + 902 p + 24 p^2 + 665 p^3 + 796 p^4 + O(p^5), \\
			x_2 &= 272 + 1462 p + 614 p^2 + 1182 p^3 + 1912 p^4 + O(p^5),\\
			&\cdots \\
			x_6 &= 651 + 750 p + 1492 p^2 + 743 p^3 + 204 p^4 + O(p^5),\\
			&\cdots  \\
			x_9 &= 1054 + 367 p + 598 p^2 + 1629 p^3 + 521 p^4 + O(p^5),\\
			&\cdots \\
			x_{12} &= 1759 + 843 p + 1836 p^2 + 900 p^3 + 768 p^4+ O(p^5).
		\end{aligned}$$
		
		Recall the polynomial $g$ we computed above, $L = \mathbb{Q}(x_1)$ and $K = \mathbb{Q}(y)$ with $y=g(x_1)$. For $\sigma\in G = \Gal(M/\mathbb{Q})$, $$\sigma\in H \iff \sigma(y) = y \iff g(\sigma(x_1)) = g(x_1).$$
		
		Now one does a very explicit computation, to see that $g(x_i) = g(x_1) + O(p^5) \iff i\in \{1,6,9\}$. Therefore for $\sigma\in G$, those in $H$ should be (if we ignore the error $O(p^5)$) characterized by $\sigma(1) \in \{1,6,9\}$.

		What remains is a pure group-theoretical computation, so we use GAP here. 
		\begin{verbatim}
			G := Group((2,12)(4,7)(8,10),(2,3,4)(5,10,7)(8,11,12),
			(1,7)(2,9)(3,4)(5,10)(6,8)(11,12),(1,5)(2,4)(3,9)(6,11)(7,10)(8,12),
			(3,5,11)(4,12,10),(2,7,8)(4,12,10),(1,6,9)(2,8,7)(3,11,5)(4,10,12));;
			H := Filtered(G, g -> OnPoints(1,g) in [1,6,9]);
		\end{verbatim}
		
		One checks $H$ is indeed a group by typing \texttt{AsGroup(H);} and verify it has the expected index $(=4)$ in $G$: \texttt{Size(G)/Size(H)}. This proves $H = \Gal(M/K)$ is indeed the group we want to find. We have all ingredients to compute $\varrho(\beta)$ from the formula
		$$\varrho(\beta) = \frac{1}{|G|} \sum_{C\subset G} |C| \left( \frac{[G:H] |H\cap C|}{|C|}\right)^{2\beta}.$$
		
		In GAP, we list all conjugacy classes of $G$ via \texttt{ccl := ConjugacyClasses(G);}, there are 17. $|C|$ for each of them is \texttt{List(ccl,Size);}, $[G:H] |H\cap C|/|C|$ for each of them is \texttt{List(ccl, c -> 4*Size(Intersection(H,c))/Size(c));}, combining gives finally our long-sought formula $$\varrho(\beta) = \frac{1}{24}(8+6\times 2^{2\beta} + 4^{2\beta})$$
		Note that this in contrast to $K/\mathbb{Q}$ Galois case, which is simply $[K:\mathbb{Q}]^{2\beta-1}$. 
	\end{example}
	
	When $K/\mathbb{Q}$ is non-Galois, it would be an interesting topic of further investigation to find an easier way to compute $\varrho(\beta)$.
	~\\[0.05in]

	\section{Cuspidal partial zeta series}
	\subsection{Criterion of being cuspdial}
	We consider the following situation: $L/K$ is an abelian extension of number fields, with both $K$ and $L$ being Galois over $\mathbb{Q}$. Let $$G = \Gal(L/\mathbb{Q}) \qquad N = \Gal(L/K) \qquad Q = \Gal(K/\mathbb{Q}).$$ That is, we have an exact sequence of groups:
	$$1\longrightarrow N\longrightarrow G\longrightarrow Q\longrightarrow 1$$
	
	It determines an action of $Q$ on $N$, hence also on $\widehat{N}$. Explicitly: for $\tau \in Q$, let $\widetilde{\tau}$ be a lift in $G$, define $\tau \chi$ as $$(\tau \chi)(h) = \chi(\widetilde{\tau}^{-1} h \widetilde{\tau}), \qquad h\in N.$$
	Evidently it's independent of the lift $\widetilde{\tau}$. For character $\chi$ on $N$, denote $\chi^\ind$ as the induced character of $\chi$ to $G$, explicitly:
	$$\chi^\ind(g) = \begin{cases} 0 \qquad & g\notin N \\ \sum_{s\in G/H} \chi(s^{-1}gs) = \sum_{\tau\in Q} (\tau\chi) (g) \qquad & g\in N \end{cases}.$$
	
	\begin{lemma}
		For $\chi_1, \chi_2 \in \widehat{N}$, $\langle \chi_1^\ind, \chi_2^\ind \rangle_G$ equals the number of $\tau \in Q$ such that $\chi_1 = \tau \chi_2$.
	\end{lemma}
	\begin{proof}We have
		$$\begin{aligned}\langle \chi_1^\ind, \chi_2^\ind \rangle_G &= \frac{1}{|G|}\sum_{h\in N}  \chi_1^\ind(h) \overline{\chi_2^\ind(h)} \\
			&= \frac{1}{|G|}\sum_{h\in N}\sum_{\tau_1,\tau_2\in Q} (\tau_1 \chi_1 \overline{\tau_2 \chi_2})(h) \\
			&= \frac{|N|}{|G|} \# \{(\tau_1,\tau_2)\in Q^2 | \tau_1 \chi_1 = \tau_2 \chi_2\}  = \frac{|N||Q|}{|G|} |\{\tau\in Q | \chi_1 = \tau \chi_2\}|
		\end{aligned},$$
		then $|G| = |N||Q|$ implies the statement.
	\end{proof}
	
	Recall our notations on the Artin $L$-functions $L(\chi,s), L(\chi^\ind,s)$ and the number $a(\chi,n)$ with $L(\chi,s) = \sum_{n\geq 1} a(\chi,n)/n^s$ from first section.

	\begin{proposition}\label{cusp_criterion_prop}
		For $\chi \in \widehat{N}$, the following are equivalent:
		\begin{enumerate}
			\item $\chi^\ind$ is an irreducible representation of $G$
			\item $\sum_{n\leq x} |a(\chi,n)|^2 \ll x$
			\item For each $\tau\neq 1$ in $Q$, $\tau \chi \neq \chi$.
		\end{enumerate}
	\end{proposition}
	\begin{proof}
		$(1)\iff (3)$ follows from above lemma. Recall that in Proposition \ref{characterorder}, we have shown
		$$\sum_{n\leq x} |a(n,\chi)|^{2\beta} \ll x (\log x)^{\varrho(\chi,\beta)-1},$$
		and since $\varrho(\chi,1)$ is exactly $\langle \chi^\ind, \chi^\ind\rangle_G$, we have $(1)\iff (2)$. 
	\end{proof}
	
	Consider the Artin $L$-function $L(\chi,s)$ coming from an automorphic representation of $GL(1)/K$, by automorphic induction, there should be an automorphic representation $\pi$ on $GL(n)/\mathbb{Q}, n = [K:\mathbb{Q}]$ realizing $L(\chi,s)$. Using the growth rate $\sum_{n\leq x} |a(n,\chi)|^2 \ll x$ above, we see that $\pi$ should be cuspidal if and only if $\chi^\ind$ is irreducible. Proving this rigorously would be very hard, and is known only when $K/\mathbb{Q}$ is cyclic \cite[Section~3.6]{arthur1989simple}\footnote{the criterion there is: $\chi \in \widehat{N}$ is non-cuspidal if and only if it is invariant under some element of the Galois group, this is our condition $(3)$ of the Proposition}. We simply take above three equivalent criteria as a \textit{working definition} of cuspidal $L$-functions for $L(\chi,s)$. 
	
	% \begin{remark}
		% One can drop the hypothesis $L/K$ being abelian, requiring it to be merely Galois. Define $\widehat{N}$ to be then the set of irreducible representations of $N$. The statement of Proposition \ref{cusp_criterion_prop} remains valid. Let $\chi$ be an irreducible representation of $N$ of dimensional $n$, so $L(\chi,s)$ corresponds to an automorphic form on $\text{GL}(n)/K$, performing induction to obtain a form on $\text{GL}(n |Q|)/\mathbb{Q}$, the Proposition should detect whether it comes from a cuspidal form on $\text{GL}(n |Q|)/\mathbb{Q}$. We will not work with the more general situation $L/K$ being non-abelian here. 
		% \end{remark}
	
	\begin{remark}
		If $L$ is Hilbert class field of $K$, then we can identify $\widehat{N}$ as character group of the ideal class group, then action of $\tau \in Q = \Gal(K/\mathbb{Q})$ is the same as action of $\tau$ on ideal classes of $K$: $(\tau \chi)(\mathfrak{A}) = \chi(\mathfrak{A}^\tau)$ for an ideal class $\mathfrak{A}$ in $K$. This is because the isomorphism between class group and $\Gal(L/K)$ is given by Frobenius map, and $\tau^{-1} (\Frob \mathfrak{p}) \tau = \Frob (\tau\mathfrak{p})$. 
	\end{remark}
	
	\begin{example}
		When $K$ is quadratic, the non-trivial element of $Q$ acts as inverse on ideal class because $1+\tau$ annihilates the class group. Therefore $L(s,\chi)$ should be a cusp form of $GL(2)/\mathbb{Q}$ if and only if $\tau\chi = \bar{\chi} \neq \chi$, i.e. if and only if $\chi$ is not real.
	\end{example}
	
	\begin{example}\label{-23_disc_example}
		Let $K=\mathbb{Q}(\sqrt{-23})$, it has class number $3$, for non-real character $\chi$ of its ideal class group, above example says $L(\chi,s)$ should come from a cuspidal automorphic form on $GL(2)/\mathbb{Q}$. \par Indeed, $L(\chi,s) = L(f,s)$ with $f$ is the unique weight $1$ normalized newform of level $23$, explicitly $f(z) = \eta(z)\eta(23z)$ with $\eta(z)$ the Dedekind eta function. This particular case can be shown without deep facts, see Zagier \cite[page~42-43]{bruinier20081} for the computations.\par 
		
		Above is an example of a celebrated result of Deligne and Serre \cite{deligne1974formes}: in terms of $L$-function, odd $2$-dimensional irreducible Galois representations and weight $1$ newforms are equivalent.
	\end{example}
	
	\begin{lemma}
		For $\tau\in Q, \chi\in \widehat{N}$, we have $L(\tau\chi,s) = L(\chi,s)$.
	\end{lemma}
	\begin{proof}
		This holds for any representation $\chi : G\to GL(V)$, not necessarily one-dimensional. The proof is purely formal. By definition $$L(\chi,s) = \prod_{\p} \frac{1}{\det(1-\chi(\Frob \mathfrak{P}/\p) N(\p)^{-s} | V^{I_{\mathfrak{P}/\p}})}.$$
		Here $V^{I_{\mathfrak{P}/\p}}$ is the subspace fixed under $\chi$ by inertia group $I_{\mathfrak{P}/\p}$. Note that $$\tau\chi (\Frob \mathfrak{P}/\p) = \chi (\Frob \mathfrak{\widetilde{\tau}P}/\widetilde{\tau}\p),$$
		and the fixed space under action of $\tau\chi$ is $V^{I_{\widetilde{\tau}\mathfrak{P}/\widetilde{\tau}\p}}$ therefore
		$$L(\tau\chi,s) = \prod_{\p} \frac{1}{\det(1-\chi(\Frob \mathfrak{\widetilde{\tau}P}/\widetilde{\tau}\p ) N(\p)^{-s} | V^{I_{\mathfrak{\widetilde{\tau}P}/\widetilde{\tau}\p}})}.$$
		As $K/\mathbb{Q}$ is Galois, $\widetilde{\tau}$ permutes prime ideals of $K$ and norm is Galois invariant, so they're equal. 
	\end{proof}
	
	Let $$\widehat{N}_0 = \{\chi\in \widehat{N} \mid \chi^\ind \text{ is irreducible} \},$$
	then $Q$ acts on $\widehat{N}_0$ without fix point (by third point of Proposition \ref{cusp_criterion_prop}). Let $\chi_1,\cdots,\chi_k$ be representatives from different orbits, we have $k = |\widehat{N}_0|/|Q|$. For $i\neq j$, Proposition \ref{cusp_criterion_prop} implies $\chi_i^\ind$ and $\chi_j^\ind$ are two non-isomorphic irreducible representations of $G$. Moreover, by above lemma, the Artin $L$-function $L(\chi_i,s)$ depends only on the orbit, not on representatives chosen. 
	
	\begin{proposition}\label{L_func_linear_indep_prop}
		Let $\chi_1,\cdots,\chi_k$ as above, then $L(\chi_1,s),\cdots,L(\chi_k,s)$ are linearly independent over $\mathbb{C}$.
	\end{proposition}
	\begin{proof}
		Let $\sum_i c_i L(\chi_i,s) = 0$, for $p$ unramified in $L$, comparing $p^{-s}$ coefficient gives $\sum_i c_i \chi_i^\ind(\Frob p) = 0$. Since the Frobenius map is surjective when $p$ varies, $\sum_i c_i \chi^\ind_i = 0$, $\chi^\ind_i$ being mutually different irreducible representations implies $c_i = 0$.\end{proof}

	\begin{example}
		Let $K$ be the cubic field in Example \ref{cubicexample1}, then we computed that $\varrho(\chi,\beta = 1) = 3$ for all $\chi \in \widehat{N}$, so none of $L(\chi,s)$ will be cuspidal. In this case, $|\widehat{N}_0| = 0$. \end{example}
	
	\begin{example}
		For $K$ be the cubic field in Example \ref{cubicexample2}, here $\widehat{N} = \mathbb{Z}/2\mathbb{Z}\times \mathbb{Z}/2\mathbb{Z}$, $$\varrho(\chi,\beta=1) =  \begin{cases} 3 & \chi = 1 \\ 1& \chi\neq 1\end{cases},$$
		so for $\chi\neq 1$, $\chi^\ind$ is irreducible. In this case $|\widehat{N}_0| = 3$ and there is only one orbit. This $L$-function, namely
		$$ 1-\frac{1}{5^s}+\frac{1}{8^s}-\frac{1}{13^s}-\frac{1}{17^s}-\frac{1}{23^s}+\frac{2}{25^s}+\frac{1}{27^s} + \cdots,$$
		should come from a cuspidal form on $GL(3)/\mathbb{Q}$.
	\end{example}

	\begin{remark} We make two remarks, they are however not used in sequel. \footnote{However, Example \ref{order21_example} below is a special case of these two general facts.}\par
		Firstly, if $1\to N\to G\to Q\to 1$ splits (i.e. $G$ is a semidirect product of $N,Q$), the representations $\chi_1^\ind,\cdots, \chi_k^\ind$ constructed above are exactly those irreducible representations of $G$ with highest dimension ($=|Q|$). This follows from a general result of linear representation of semidirect product (Serre, \cite[Chapter~8.2]{serre1977linear}). \par
		When $L$ is Hilbert class field of $K$ and $K/\mathbb{Q}$ cyclic, global class field theory says above sequence always splits. So in this case, $G$ as an abstract group alone determines $\varrho_\text{cusp}(\sigma,\beta)$ (defined in next section) that appears in asymptotic of moments. 
	\end{remark}
	
	\subsection{Cuspidal version of $a(\sigma,n)$}
	Recall the partial zeta function associated for $\sigma\in N$:
	$$\zeta(\sigma,s) = \frac{1}{|N|}\sum_{\chi \in \widehat{N}} \bar{\chi}(\sigma) L(\chi,s)= \sum_{n\geq 1} \frac{a(\sigma,n)}{n^s}.$$
	We remove those terms which are not cuspidal $L$-function, obtaining
	$$\zeta_\text{cusp}(\sigma,s) = \frac{1}{|N|}\sum_{\chi \in \widehat{N}_0} \bar{\chi}(\sigma) L(\chi,s) := \sum_{n\geq 1} \frac{a_\text{cusp}(\sigma,n)}{n^s}.$$
	Using the fact that $L(\chi,s)$ is invariant under action of $Q$ and $\chi^\ind = \sum_{\tau \in Q} \tau\chi$, one easily sees
	$$\zeta_\text{cusp}(\sigma,s) = \frac{1}{|N|}\sum_{i=1}^k \overline{\chi_i^\ind(\sigma)} L(\chi_i,s),$$
	here $\chi_1,\cdots,\chi_k$ are representatives of $Q$-orbits on $\widehat{N}_0$. This form turns out to be more amendable than original definition. \par
	
	We wish to investigate the moment $\sum_{n\leq x} |a_{\text{cusp}}(\sigma,n)|^{2\beta}$. Before that, we caution there might exist $\sigma \in N$ such that $\zeta_\text{cusp}(\sigma,s)$ is identically zero. By linear independence (Proposition \ref{L_func_linear_indep_prop}) of $L(\chi_i,s)$, this occurs if and only if $\chi_i^\ind(\sigma)=0$ for all $i$.
	
	\begin{example}\label{non_vanishing_form_example}
		Assume $K/\mathbb{Q}$ be quadratic, if the class group $N$ is $2$-torsion, then $\widehat{N}_0$ is empty, so $\zeta_\text{cusp}(\sigma,s) = 0$ for all $\sigma\in N$. \par
		
		Quite unexpectedly, there are examples for non-$2$-torsion $N$ and $\sigma\in N$ such that $\zeta_\text{cusp}(\sigma,s) = 0$. To be precise, note that
		$\chi_i^\ind = \chi_i + \bar{\chi}_i = 2\Re(\chi_i)$, and any $\chi_i \in \widehat{N}_0$ if and only if its order is not $1$ or $2$, so if $\sigma$ is killed by all $\chi_i^\ind$, then it must have order $4$. From which it follows that $\zeta_\text{cusp}(\sigma,s) = 0$ if and only if $N \cong \mathbb{Z}/4\mathbb{Z} \times (\mathbb{Z}/2\mathbb{Z})^n$ and $\sigma$ an element of order $4$ in $N$. 
	\end{example}
	
	Consequently, unlike the situation for $\sum_{n\leq x} |a(\sigma,n)|^{2\beta}$, whose leading term of asymptotic is independent of $\sigma \in N$, the leading term of $\sum_{n\leq x} |a_\text{cusp}(\sigma,n)|^{2\beta}$ might depend on $\sigma\in N$. Let $$\varrho_\text{cusp}(\sigma,\beta) = \max_{\substack{i=1,\cdots,k \\ \chi_i^\ind(\sigma)\neq 0}} \varrho(\chi_i,\beta),$$
	here we used the notation $$\varrho(\chi_i,\beta) = \frac{1}{|G|}\sum_{g\in G} |\chi_i^\ind(g)|^{2\beta},$$
	which is the exponent for $\sum_{n\geq 1} |a(\chi,n)|^{2\beta}/n^s$ at $s=1$.
	
	Also recall the notation $$\varrho(\chi_1,\cdots,\chi_i) = \frac{1}{|G|}\sum_{g\in G} \chi_1^\ind(g)\cdots \chi_i^\ind(g),$$
	which is always a non-negative integer and is the exponent for $\sum_{n\geq 1} a(\chi_1,n)\cdots a(\chi_i,n)/n^s$ at $s=1$. Note in both cases, the exponent at $s=1$ remaimight dependns unchanged if we only sum over squarefree integer $n$, this follows immediately from their proofs Proposition \ref{characterorder}.
	
	\begin{theorem}\label{cusporderintegral}
		Let $\beta$ be a positive integer, for $\sigma \in N$, if $\zeta_\text{cusp}(\sigma,s)$ is not identically zero, then $$\sum_{n\geq 1}\frac{|a_\text{cusp}(\sigma,n)|^{2\beta}}{n^s} $$ has pole of exact order $\varrho_\text{cusp}(\sigma,\beta)$ at $s=1$. 
		Hence $$\sum_{n\leq x}|a_\text{cusp}(\sigma,n)|^{2\beta} \sim C x (\log x)^{\varrho_\text{cusp}(\sigma,\beta) -1}$$
		for some positive constant $C$.
	\end{theorem}
	
	The proof for $\beta = 1$ is especially elegant: in this case, we need to show  $\sum_{n\geq 1} \frac{|a_{\text{cusp}}(\sigma,n)|^2}{n^s}$ has pole of exact order $1$, 
	$$\sum_{n\geq 1} \frac{|a_{\text{cusp}}(\sigma,n)|^2}{n^s} = \frac{1}{|N|^2}\sum_{i,j} \overline{\chi_i^{\ind}(\sigma)} \chi_j^\ind(\sigma) \sum_{n\geq 1} \frac{a(\chi_i,n) a(\overline{\chi_j},n)}{n^s}.$$
	For the indices $i\neq j$, the order of poles of the sum is $\varrho(\chi_i,\overline{\chi_j}) = \langle \chi_i^\ind, \chi_j^\ind \rangle_G$ which is $0$ since $\chi_i^\ind$ are non-isomorphic representations. So the order of poles of $\sum_{n\geq 1} \frac{|a_{\text{cusp}}(\sigma,n)|^2}{n^s}$ is the same as
	$$\frac{1}{|N|^2}\sum_{i} |\chi_i^{\ind}(\sigma)|^2  \sum_{n\geq 1} \frac{|a(\chi_i,n)|^2}{n^s},$$
	each has pole order $\langle \chi_i^\ind, \chi_i^\ind \rangle_G = 1$, and by assumption the projection onto cuspidal space is non-zero $\iff$ at least one $|\chi_i^\ind(\sigma)| > 0$, so the leading term is positive, completing the proof when $\beta = 1$.
	
	For general positive integral $\beta$, we need some notations and easy observations. Fix $\sigma \in N$, let $\mathfrak{X}_{\sigma,\beta}$ be those $\chi_i$ such that maximum is attained for $\varrho_{\text{cusp}}(\sigma,\beta)$, also let $\widehat{N}_{0,\sigma}$ be those $\chi_i$ in $\widehat{N}_0$ for which $\chi_i^\ind(\sigma)\neq 0$. 
	\begin{lemma}
		Let $\{\rho_1,\cdots,\rho_r\}, \{\psi_1,\cdots,\psi_r\} \subset \widehat{N}_0$ such that $$\varrho(\rho_1,\cdots,\rho_r,\overline{\psi_1},\cdots,\overline{\psi_r}) = \varrho_{\text{cusp}}(\sigma,r).$$
		Then
		\begin{enumerate}
			\item each $\rho_i,\psi_i \in \mathfrak{X}_{\sigma,r},$
			\item $|\rho_i^\ind(g)| = |\psi_j^\ind(g)|$ for all $g\in G,$
			\item $\rho^\ind_1 \cdots \rho^\ind_r = \psi^\ind_1 \cdots \psi^\ind_r.$
		\end{enumerate}
	\end{lemma}
	\begin{proof}Hölder's inequality implies
		$$\begin{aligned}\varrho(\rho_1,\cdots,\rho_r,\overline{\psi_1},\cdots,\overline{\psi_r}) &= \frac{1}{|G|} \sum_{g\in G} \rho^\ind_1(g) \cdots \rho^\ind_r(g)\overline{\psi^\ind_1}(g)\cdots \overline{\psi^\ind_r}(g) \\ &\leq \frac{1}{|G|} \sum_{g\in G} |\rho^\ind_1(g)| \cdots |\rho^\ind_r(g)| |\psi^\ind_1(g)|\cdots |\psi_r^\ind(g)|  
			\\ &\leq \left(\frac{1}{|G|} \sum_{g\in G} |\rho_1^\ind(g)|^{2r} \right)^{1/2r}\cdots \left(\frac{1}{|G|} \sum_{g\in G} |\psi_r^\ind(g)|^{2r} \right)^{1/2r} = \varrho(\rho_1,r)^{1/{2r}} \cdots \varrho(\psi_r,r)^{1/2r}
		\end{aligned}.$$
		In order for this to equal to $$\varrho_\text{cusp}(\sigma,\beta) = \max_{\substack{i=1,\cdots,k \\ \chi_i^\ind(\sigma)\neq 0}} \varrho(\chi_i,\beta),$$
		each of term inside the parenthesis must be this number, proving $(1)$. During the Hölder's inequality step, equality holds if and only if there exists $\lambda$ such that $\forall g\in G, |\phi_1^\ind(g)| = \lambda|\phi_2^\ind(g)|$, here $\phi_1,\phi_2 \in \{\rho_1,\cdots,\rho_r,\psi_1,\cdots,\psi_r\}$, because $\phi_1,\phi_2$ are both $1$-dim characters, $\phi_1^\ind(1) = \phi_2^\ind(1) = |Q|$, thus $\lambda = 1$, giving $(2)$. \par
		Using polar coordinates, write $$(\rho^\ind_1 \cdots \rho^\ind_r)(g) = R(g) e^{i \theta_1(g)}, \qquad (\psi^\ind_1 \cdots \psi^\ind_r)(g) = R(g) e^{i \theta_2(g)},$$
		then $$\frac{1}{|G|} \sum_{g\in G} R(g)^2 e^{i\theta_1(g) - i\theta_2(g)} = \varrho(\rho_1,\cdots,\rho_r,\overline{\psi_1},\cdots,\overline{\psi_r}) = \varrho_{\text{cusp}}(\sigma,r) = \frac{1}{|G|}\sum_{g\in G} R(g)^2, $$
		this forces $e^{i\theta_1(g) - i\theta_2(g)} = 1$, which is $(3)$. 
	\end{proof}
	
	\begin{proof}[Proof of Theorem \ref{cusporderintegral}]
		Because $a_\text{cusp}(\sigma,n)$ is a linear combination of $a(\chi_i,n)$, it is easy to see order of pole of both Dirichlet series are $\leq \varrho_\text{cusp}(\sigma,\beta)$. Therefore it suffices to prove the pole for squarefree sum $$\sum_{n\geq 1, n \text{ squarefree}}\frac{|a_\text{cusp}(\sigma,n)|^{2\beta}}{n^s}$$ is exactly this number. Let $$X = \{(\rho_1,\cdots,\rho_\beta, \psi_1,\cdots,\psi_\beta) \in (\widehat{N}_{0,\sigma})^{2\beta} \mid \varrho(\rho_1,\cdots,\rho_\beta,\overline{\psi_1},\cdots,\overline{\psi_\beta}) = \varrho_{\text{cusp}}(\sigma,\beta)\}$$
		be the set of $2\beta$-tuples that attain the maximum order of pole, this set is non-empty since we assumed $a_\text{cusp}(\sigma,n)$ does not vanish identically. Furthermore,
		\begin{multline*}\sum_{n\geq 1} \frac{|a_{\text{cusp}}(\sigma,n)|^{2\beta}}{n^s} = \frac{1}{|N|^{2\beta}}\sum_{(\rho_i,\psi_i)\in (\widehat{N}_{0,\sigma})^{2\beta}} (\rho^\ind_1 \cdots \rho^\ind_\beta)(\sigma)\overline{(\psi^\ind_1 \cdots \psi^\ind_\beta)(\sigma)} \\ \times \sum_{n\geq 1} \frac{\prod_{1\leq i\leq \beta} a(\rho_i,n) \prod_{1\leq i\leq \beta} a(\overline{\psi_i},n)}{n^s}.\end{multline*}
		Consider whether $(\rho_i,\psi_i)$ is in $X$ or not, when no, it only contributes pole of lower order, so we can safely ignore them, for those in $X$. By lemma, we have $\rho^\ind_1 \cdots \rho^\ind_\beta = \psi^\ind_1 \cdots \psi^\ind_\beta$, so modulo terms of lower order poles. Thus
		$$\sum_{n\geq 1} \frac{|a_{\text{cusp}}(\sigma,n)|^{2\beta}}{n^s} \sim \frac{1}{|N|^{2\beta}}\sum_{(\rho_i,\psi_i)\in X} |(\rho^\ind_1 \cdots \rho^\ind_\beta)(\sigma)|^2 \sum_{n\geq 1} \frac{\prod_{1\leq i\leq \beta} a(\rho_i,n) \prod_{1\leq i\leq \beta} a(\overline{\psi_i},n)}{n^s}.$$
		The same is true if we restrict the sum to square-free $n$ only, in this case, the $a(\chi_i,n)$ becomes multiplicative:
		$$\prod_{1\leq i\leq \beta} a(\rho_i,n) \prod_{1\leq i\leq \beta} a(\overline{\psi_i},n) = a(\rho^\ind_1\otimes \cdots \otimes \rho^\ind_\beta,n) \overline{a(\psi^\ind_1\otimes \cdots \otimes \psi^\ind_\beta,n)}= |a(\rho^\ind_1\otimes \cdots \otimes \rho^\ind_\beta,n)|^2,$$
		which are coefficient associated with tensor product representation. Therefore all terms in above displayed equation have order of pole $\varrho_{\text{cusp}}(\sigma,r)$ and positive leading coefficients, therefore so is $\sum_{n\geq 1} \frac{|a_{\text{cusp}}(\sigma,n)|^{2\beta}}{n^s}$.
	\end{proof}
	
	\begin{remark}
		The above proof can be directly adapted to prove Theorem \ref{partialorder} when $\beta$ is positive integral.
	\end{remark}
	
	The above proof only works for $\beta$ integral, for general real $\beta>0$, we requires an additional assumption for given $\sigma \in N, \beta>0$:
	\begin{equation}\label{characteres}\tag{**}\chi_i \in \mathfrak{X}_{\sigma,\beta} \implies \forall g\in G, \chi_i^\ind(g) \in \mathbb{R}\times \text{(root of unity)}. \end{equation}
	Hererecall that $\mathfrak{X}_{\sigma,\beta}$ is the set of character for which maximum $\varrho_\text{cusp}(\sigma,\beta) = \max_{\substack{i=1,\cdots,k \\ \chi_i^\ind(\sigma)\neq 0}} \varrho(\chi_i,\beta)$ is attained. 
	
	\begin{theorem}\label{cusporderreal}
		Let $\beta > 0$, $\sigma \in N$, if $\zeta_\text{cusp}(\sigma,s)$ is not identically zero and equation (\ref{characteres}) is satisfied, writing $\varrho = \varrho_\text{cusp}(\sigma,\beta)$, then $$0 < \liminf_{s\to 1^{-}} (s-1)^{\varrho} \sum_{n\geq 1}\frac{|a_\text{cusp}(\sigma,n)|^{2\beta}}{n^s} \leq \limsup_{s\to 1^{-}} (s-1)^{\varrho} \sum_{n\geq 1}\frac{|a_\text{cusp}(\sigma,n)|^{2\beta}}{n^s} < \infty .$$ 
	\end{theorem}
	\begin{proof}The proof goes almost exactly the same as that of Theorem \ref{partialorder}.
		
		The fact that $\limsup < \infty$ is evident. Proving $\liminf \neq 0$ is more involved. We start with $$\zeta_\text{cusp}(\sigma,s) = \frac{1}{|N|}\sum_{i=1}^k \overline{\chi_i^\ind(\sigma)} L(\chi_i,s).$$
		Abbreviate $\chi^\ind(C) := a_{C,\chi}$ where $C$ is a conjugacy class of $G$. By our assumption (\ref{characteres}), we can choose positive integer $l$ such that $a_{C,\chi}^l \geq 0$ for all $C$ and $\chi \in \mathfrak{X}_{\sigma,\beta}$, write $\mu$ as a primitive $l$-th root of unity. For some $\rho_i(C)\in \mathbb{C}, i=0,\cdots,l-1$ to be fixed later and $\chi \in \{\chi_1,\cdots,\chi_k\}$, define
		$$f(\chi,C) = \sum_{i=0}^{l-1} \rho_i(C) \prod_{\Frob p \in C} (1+ \frac{\mu^i a_{C,\chi}}{p^s}).$$ 
		Also define $a'(\chi,n)$ and $a'(\sigma,n)$ via
		$$\sum_{n\geq 1} \frac{a'(\chi,n)}{n^s} := \prod_{C\subset G} f(\chi,C),$$
		$$a'(\sigma,n) = \frac{1}{|N|}\sum_{i=1}^k \overline{\chi_i^\ind(\sigma)} a'(\chi_i,n).$$
		Note that the same prime $p$ never occurs in two $f(\chi,C_1)$ and $f(\chi,C_2)$ with $C_1\neq C_2$; as usual, we only focus on prime which are not ramified in $L$.
		\\[0.2in]
		We claim that we can choose $\rho_i(C) \neq 0$ independent of $\chi$ such that 
		\begin{enumerate}
			\item For all $j$ and $C$, $\sum_{i=0}^{l-1} \rho_i(C) \mu^{ij}$ equals $0$ or $1$.
			\item $\overline{\chi_i^\ind(\sigma)}a'(\chi_i,n) \geq 0$ for all $n$ and $i$.
		\end{enumerate}
		
		Assuming above claim, let us prove $\liminf > 0$. Since $p$-th coefficient of $L(\chi_i,s)$ is $\chi_i^\ind(\Frob p)$, (1) implies $a'(\chi,n) = a_\text{cusp}(\chi,n)$ for all $\chi$ or $a'(\chi,n) = 0$ for all $\chi$, therefore $|a_\text{cusp}(\sigma,n)| \geq |a'(\sigma,n)|$, so it suffices to prove the assertion on $\liminf$ for $\sum_{n\geq 1} |a'(\sigma,n)|^{2\beta}n^{-s}$. Observe
		$$\begin{aligned}|a'(\sigma,n)|^{2\beta} &\gg \left|\sum_{\chi \in \mathfrak{X}_{\sigma,\beta}} \overline{\chi_i^\ind(\sigma)} a'(\chi,n) \right|^{2\beta} - \sum_{\chi \notin \mathfrak{X}_{\sigma,\beta}} |a'(\chi,n)|^{2\beta} \\
			&\gg \sum_{\chi \in \mathfrak{X}_{\sigma,\beta}} |a'(\chi,n)|^{2\beta} - \sum_{\chi \notin \mathfrak{X}_{\sigma,\beta}} |a'(\chi,n)|^{2\beta} \qquad \text{by condition (2)}.\end{aligned}$$
		Condition (1) also implies (via the same reasoning as in Theorem \ref{partialorder})
		$$\sum_{n\geq 1} \frac{|a'(\chi,n)|^{2\beta}}{n^s} = \prod_{C\subset G} \left( \sum_{i=0}^{l-1} \rho_i(C) \prod_{\Frob p \in C} (1+ \frac{\mu^i |a_{C,\chi}|^{2\beta}}{p^s}) \right).$$
		each $\prod_{\Frob p \in C} (1+ \frac{\mu^i |a_{C,\chi}|^{2\beta}}{p^s})$ is $\asymp (s-1)^{-\mu^i |a_{C,\chi}|^{2\beta} |C|/|G|}$ as $s\to 1$, so as $i$ varies between $0$ and $l-1$, the dominant term will be $\asymp (s-1)^{-|a_{C,\chi}|^{2\beta} |C|/|G|}$ (here we also used $\rho_i(C)\neq 0$), multiplying over all $C\subset G$ implies that RHS of above displayed equation is $\asymp (s-1)^{-\varrho(\chi,\beta)}$, which is  $(s-1)^{-\varrho_\text{cusp}(\sigma,\beta)}$ when $\chi \in \mathfrak{X}_{\sigma,\beta}$, 
		therefore $$\sum_{n\geq 1} \frac{|a'(\sigma,n)|^{2\beta}}{n^s} \gg (s-1)^{-\varrho_\text{cusp}(\sigma,\beta)} - \sum_{\chi \notin \mathfrak{X}_{\sigma,\beta}} \sum_{n\geq 1} \frac{|a'(\chi,n)|^{2\beta}}{n^s}.$$
		By definition of $\mathfrak{X}_{\sigma,\beta}$, all terms on the right has pole order $<\varrho_\text{cusp}(\sigma,\beta)$, so above is still $\gg (s-1)^{-\varrho_\text{cusp}(\sigma,\beta)}$, proving $\liminf$ is positive, assuming two conditions above. 
		\\[0.2in]
		Now we explain how to choose $\rho_i(C)$ achieving the criteria. We make
		\begin{equation}\label{rhochoice2}\rho_i(C) = l^{-1} \quad \forall i \qquad   \qquad\text{ or } \qquad \rho_i(C) = l^{-1} \mu^{-1-i} \quad \forall i.\end{equation}
		Obviously $(1)$ is satisfied. For $(2)$, let $C_0$ be the conjugacy class containing $\sigma$. We make the first choice if $\sigma \notin C$ and otherwise the second choice. For $\sigma\notin C$, 
		$$f(\chi,C) = 1 + \sum_{p_i\in \Frob C, p_i\neq p_j} \frac{a_{C,\chi}^l}{p_1^s \cdots p_l^s} + \sum_{p_i\in \Frob C, p_i\neq p_j} \frac{a_{C,\chi}^{2l}}{p_1^s \cdots p_{2l}^s} + \cdots  \qquad C\neq C_0,$$
		which has non-negative coefficients (by our choice of $l$); for $\sigma\in C$, giving
		$$f(\chi,C_0) = \sum_{p_1\in \Frob C_0} \frac{a_{C,\chi}}{p_1^s} + \sum_{p_i\in \Frob C_0, p_i\neq p_j} \frac{a_{C_0,\chi}^{l+1}}{p_1^s \cdots p_{l+1}^s} + \sum_{p_i\in \Frob C_0, p_i\neq p_j} \frac{a_{C_0,\chi}^{2l+1}}{p_1^s \cdots p_{2l+1}^s} + \cdots,$$
		so $\overline{\chi_i^\ind(\sigma)} f(\chi,C_0)$ also has non-negative coefficients. Then
		$$\sum_{n\geq 1} \frac{\overline{\chi_i^\ind(\sigma)} a'(\chi,n)}{n^s} = \overline{\chi_i^\ind(\sigma)} f(\chi,C_0) \prod_{C\neq C_0} f(\chi,C)$$ also have non-negative coefficients, which is $(2)$. 
	\end{proof}
	
	\begin{remark}
		Just as proof of Theorem \ref{partialorder}, both Theorem \ref{cusporderintegral} and \ref{cusporderreal} generalize when one replaces $\sum_{n\geq 1} \frac{|a_\text{cusp}(\sigma,n)|^{2\beta}}{n^s}$ with $$\sum_{\substack{(n,S)=1 \\ n \text{ squarefree}}} \frac{|a_\text{cusp}(\sigma,n)|^{2\beta}}{n^s}.$$
	\end{remark}

	One drawback of Theorem \ref{cusporderreal} is that (\ref{characteres}) is in general \textit{difficult to check} for a given field extension $L/K$. Even when $K/\mathbb{Q}$ is cyclic (in which automorphic induction has been proven)  and $L$ being Hilbert class field of $K$, it might still fail. 
	
	\begin{example}\label{order21_example}
		Let $p \equiv 1 \pmod{3}$ be a prime, $K$ be the unique cubic subfield of $\mathbb{Q}(\zeta_p)$. Assume $K$ has class number $7$,\footnote{this is the case when $p=313,877,1129,\cdots$} $L$ its Hilbert class field, we have an exact sequence
		$$1\longrightarrow \Gal(L/K) := N \longrightarrow \Gal(L/\mathbb{Q}) := G \longrightarrow \Gal(K/\mathbb{Q}) := Q \longrightarrow 1.$$
		$G$ is a group of order $21$, there are only two such groups up-to isomorphism: $C_7\times C_3$ or the non-abelian semidirect product $C_7 \rtimes C_3$. \par We claim $G$ cannot be $C_7\times C_3$. Indeed, if it were, then every subfield of $L$ is Galois over $\mathbb{Q}$ since $G$ would be abelian, in particular, this holds for the inertia field $F$ of $\mathfrak{P}/p$, here $\mathfrak{P}$ is a prime in $L$ lying over $p$, $F/\mathbb{Q}$ is degree $7$, $(\mathfrak{P} \cap F)/p$ is unramified in $F$, since $F$ would also be Galois, it must be unramified for every prime lying above $p$, so $F/\mathbb{Q}$ is everywhere unramified, a contradiction to Minkowski's theorem. \par
		Thus $G$ is $C_7 \rtimes C_3$, it has two $3$-dimensional characters:
		
		\begin{table}[H]
			\centering
			\begin{tabular}{@{}cccccc@{}}
				\toprule
				$\text{conjugacy class (size)}$ & (1) & (7) & (7) & (3) & (3) \\ \midrule
				$\rho_1$ & 3 & 0 & 0 & $\frac{-1-\sqrt{-7}}{2}$ & $\frac{-1+\sqrt{-7}}{2}$ \\
				$\rho_2$ & 3 & 0 & 0 & $\frac{-1+\sqrt{-7}}{2}$ & $\frac{-1-\sqrt{-7}}{2}$ \\ \bottomrule
			\end{tabular}
			\caption{The two $3$-dimensional characters of $G$, the ideal class group $\Gal(L/K)$ is union of 1st, 4th and 5th conjugacy classes}
			\label{tab:my-table}
		\end{table}
		Action of $Q$ on $\widehat{N}_0 =\{ \chi\in \widehat{N} | \chi^\ind \text{ irreducible}\}$ has $2$ orbits, $\rho_i = \chi_i^\ind$, their $L$-functions $L(\rho_1,s)$ and $L(\rho_2,s)$ differ by complex conjugation. \par
		For any ideal class $\sigma$ of $N$, we have $$\sum_{n\geq 1} \frac{a_\text{cusp}(\sigma,n)}{n^s} = \frac{1}{7} \sum_{i=1,2} \overline{\rho_i(\sigma)} L(\rho_i,s).$$ Looking at the table, we see none of the cuspidal projection vanishes. When $\beta\in \mathbb{Z}^{\geq 1}$, each of them has growth:
		$$\sum_{n\leq x} |a_\text{cusp}(\sigma,n)|^{2\beta} \sim C x (\log x)^{\varrho_\text{cusp}(\sigma,\beta)-1},$$
		where $$\varrho_\text{cusp}(\sigma,\beta) = \frac{1}{21}(3^{2\beta} + 3\left| \frac{-1-\sqrt{-7}}{2}\right|^{2\beta} + 3\left| \frac{-1+\sqrt{-7}}{2}\right|^{2\beta}) = \frac{1}{7}(3^{2\beta-1} + 2^{1+\beta}).$$
		However, for non-integral $\beta > 0$, since the prerequisite (\ref{characteres}) no longer holds: $\frac{-1\pm \sqrt{-7}}{2}$ has no desired property, so we cannot conclude from Theorem \ref{cusporderreal} the same order of growth for such $\beta$.
	\end{example}
	
	Similar to non-cuspidal case, we expect the following to hold:
	\begin{conjecture}
		For real $\beta>0$, Theorem \ref{cusporderreal} still holds after removing the condition in equation (\ref{characteres}); $\liminf$ and $\limsup$ in Theorem \ref{cusporderreal} should also be equal. 
	\end{conjecture}
	Note that, unlike the non-cuspidal case, the $\limsup$ and $\liminf$ might depend on $\sigma$.

	\subsection{Epistein zeta function of non-fundamental discriminant}
	
	In all our previous explicit examples, we focused on the case $L$ being the Hilbert class field of $K$. The theorems however works for any abelian extension of number field $L/K$. In this subsection, we shall apply our previous results to certain abelian extensions closely related to non-maximal orders of $K$. \par
	
	We recall some basic facts relating binary quadratic form to class field theory, (c.f. \cite{cox2022primes}). Let $K$ be an imaginary quadratic field, $\mathcal{O}$ an order of $K$, discriminant of $\mathcal{O}$ be $D$, $f$ the conductor of $\mathcal{O}$, $L$ be its associated ring class field (it is characterized by $p\nmid D$ splits in $L$ if and only if $p$ is represented by the principal form). $S'$ be primes of $K$ lying above $f$, let $I^{S'}$ be free abelian group generated by prime ideals in $K$ coprime to $f$. Similar to the fundamental-discriminant case, $L/\mathbb{Q}$ is still Galois. \par

	Recall the number $a(\sigma,n)$ defined previously, by examining the proof of Lemma \ref{partial_zeta_interpretation_lemma}, one sees when $(n,f)=1$, $a(\sigma,n)$ is the number of ideals in $I^{S'}$ of norm $n$ that maps to $\sigma$ under Frobenius map. 
	
	Since $G = \Gal(L/\mathbb{Q})$ has a normal subgroup $N = \Gal(L/K)$, the argument in Example \ref{quad_example} (which was done there only when $\mathcal{O} = \mathcal{O}_K$) carries through, and we have the same formula:
	$$\varrho(\chi,\beta) = \frac{1}{2m} \sum_{j=0}^{m-1} \left|2\cos(\frac{2\pi j}{m})\right|^{2\beta}, \qquad m \text{ = order of }\chi.$$
	So $\varrho(\beta) = [K:\mathbb{Q}]^{2\beta-1} = 2^{2\beta-1}$. The condition that $\chi^\ind$ is irreducible remains unchanged: it is irreducible if and only if $\chi$ is non-real. 
	\\[0.02in]
	
	$\Gal(L/K)$ is isomorphic to the \textit{form class group}: primitive quadratic form of discriminant $D$ upto equivalence. Let $g(x,y)$ be a primitive binary quadratic form of discriminant $D$, $$r_g(n) = \# \{(x,y)\in \mathbb{Z}^2 \mid g(x,y) = n\}.$$ The quantity $a(\sigma,n)$ for $(n,S)=1$ is essentially the number of $r_g(n)$ with $g$ corresponds to $\sigma \in \Gal(L/K)$. More precisely, we can choose a representative $g(x,y) = ax^2+bxy+cy^2$ such that $\mathfrak{a} = \mathbb{Z}a + \mathbb{Z}\frac{-b+\sqrt{D}}{2}$ is an ideal of $\mathcal{O}$ prime to the conductor $f$, and we have
	$$g(x,y) = \frac{N(xa + y\frac{-b+\sqrt{D}}{2})}{N(\mathfrak{a})}.$$
	So
	$$\sum_{n\geq 1}\frac{r_g(n)}{n^s} = N(\mathfrak{a})^s \sum_{z\in \mathfrak{a}} \frac{1}{N(z)^s} = N(\mathfrak{a})^s w\sum_{(z)\subset \mathfrak{a}} \frac{1}{N(z)^s},$$
	with $w$ the number of units in $\mathcal{O}$. Because ideals in $\mathcal{O}$ lacks unique factorization, it will be more convenient to look at those prime to conductor, which do have unique factorization, this part corresponds to 
	$$N(\mathfrak{a})^s\sum_{(z)\subset \mathfrak{a}, (z,f)=1} \frac{1}{N(z)^s} = \frac{1}{h(\mathcal{O})} \sum_\chi \chi(\mathfrak{a}) \sum_{I\subset \mathcal{O}, (I,f)=1} \frac{\chi(I)}{N(I)^s}.$$
	Here $\chi$ ranges over all characters of the Picard group of $\mathcal{O}$: i.e. ideal class that are prime to $f$, this is isomorphic to $\Gal(L/K)$. \par
	
	By lifting $I\subset \mathcal{O}$ to $\mathcal{O}_K$, the sum $$\sum_{I\subset \mathcal{O}, (I,f)=1} \frac{\chi(I)}{N(I)^s} = \sum_{I\subset \mathcal{O}_K, (I,f)=1} \frac{\chi(I)}{N(I)^s}$$ is the same (apart from a finite Euler product supported on $f$) as $L(\chi,s)$, with $\chi$ is now interpreted as a character on $I^{S'}$ modulo an equivalence relation: $I' \sim I$ if and only if $I'I^{-1}$ is a principal ideal of the form $(a), a\in \mathbb{Z},(a,f)=1$. \par
	
	Combining above equations, we arrive at
	$$\sum_{n\geq 1, (n,f)=1}\frac{r_g(n)}{n^s} = \frac{w}{h(\mathcal{O})} \sum_\chi \chi(\mathfrak{a}) \times (\text{some Euler factors supported at } f)\times L(\chi,s).$$
	
	We define $\sum_{n\geq 1, (n,f)=1} r_\text{cusp,g}(n) n^{-s}$ to be RHS with real $\chi$ removed, it is the $L$-function of an elliptic cusp form. 
	~\\[0.02in]
	
	Corollary \ref{ideal_counting_moment_corollary} in this case becomes:
	
	\begin{corollary}
		For real $\beta>0$, 
		$$\sum_{n\leq x, (n,f)=1} r_g(n)^{2\beta} \asymp x (\log x)^{2^{2\beta-1}-1}.$$
		With $\asymp$ replaced by $\sim$ times a constant if $\beta \in \mathbb{Z}^{\geq 1}$. 
	\end{corollary}
	
	This result looks less impressive, but it generalizes almost verbatim to the cuspidal version $r_{\text{cusp},g}(n)$: $\chi^\ind$ is always real-valued, so the prerequisite (\ref{characteres}) of Theorem \ref{cusporderreal} is satisfied. Recall the expected order of cuspidal moment: 
	$$\varrho_\text{cusp}(\sigma,\beta) = \max_{\substack{i=1,\cdots,k \\ \chi_i^\ind(\sigma)\neq 0}} \varrho(\chi_i,\beta).$$
	Here $\chi_1,\cdots,\chi_k$ are orbit representatives of $Q$ on $\widehat{N}_0$: in this case, just pick one from each pair $(\chi,\bar{\chi})$ of non-real $\chi$.
	
	\begin{corollary}[Cuspidal projection with non-fundamental discriminant]
		Assume $r_{\text{cusp},g}(n)$ is not identically zero\footnote{when this happens is described precisely in Example \ref{non_vanishing_form_example}}, let $$\rho = \varrho_\text{cusp}(\sigma,\beta) = \max_{\substack{i=1,\cdots,k \\ \chi_i^\ind(\sigma)\neq 0}} \varrho(\chi_i,\beta).$$ For any real $\beta >0$, we have
		$$\sum_{n\leq x, (n,f)=1} |r_{\text{cusp},g}(n)|^{2\beta} \asymp x (\log x)^{\rho-1}.$$
		With $\asymp$ replaced by $\sim$ times a constant if $\beta \in \mathbb{Z}^{\geq 1}$. 
	\end{corollary}
	
	The case when $D$ is a fundamental discriminant is proved in  \cite{blomer2004cusp} by using an explicit calculation. In our approach, we are able to get away with such calculations using language of finite group representation; and establish the result for all $D$, whether fundamental or not, in a uniform way. 
	~\\[0.05in]

	\section*{Acknowledgements}
	\begin{spacing}{1.3}
		I am grateful to my supervisor Prof. Valentin Blomer, for introducing me to these topics and his numerous suggestions that improve the readability of this manuscript. I also thank Dr. Egdar Assing on grading the thesis from which this manuscript originates, as well as colleague Francisco Araújo for related discussions.\\[0.1in]
	\end{spacing}

	\bibliographystyle{plain} % We choose the "plain" reference style
	\bibliography{ref.bib} % Entries are in the refs.bib file
	
\end{document}